\documentclass[reqno,11pt]{amsart}
\usepackage{amsmath,amsthm}
\usepackage{latexsym}
\usepackage{amssymb}
\usepackage[all,cmtip]{xy}
\newtheorem{Proposition}{Proposition}[section]
\newtheorem{Definition}[Proposition]{Definition}
\newtheorem{Lemma}[Proposition]{Lemma}
\newtheorem{Theorem}[Proposition]{Theorem}
\newtheorem{MainTheorem}{Theorem}
\newtheorem{Corollary}[Proposition]{Corollary}

\DeclareMathOperator{\p}{\mathbb{P}}
\newcommand{\R}{\mathbb{R}}
\newcommand{\C}{\mathbb{C}}

\DeclareMathOperator{\spt}{spt}

\DeclareMathOperator{\WF}{WF}
\DeclareMathOperator{\Ker}{Ker}
\DeclareMathOperator{\supp}{supp}
\DeclareMathOperator{\Val}{Val}
\DeclareMathOperator{\Dens}{Dens}
\DeclareMathOperator{\vol}{vol}
\DeclareMathOperator{\im}{im}
\DeclareMathOperator{\grad}{grad}

\newcommand\delbundle[1]{T^*#1 \setminus \underline{0}}
\usepackage{srcltx}
\usepackage{xcolor}

\title{Convolution of valuations on manifolds}
\author{Semyon Alesker}
\author{Andreas Bernig}

\email{semyon@post.tau.ac.il}
\email{bernig@math.uni-frankfurt.de}

\address{School of Mathematical Sciences, Tel Aviv University, 69978 Tel Aviv, Israel}
\address{Institut f\"ur Mathematik, Goethe-Universit\"at Frankfurt, Robert-Mayer-Str. 10, 60054
Frankfurt, Germany}
\begin{document}

\begin{abstract}
We introduce the new notion of convolution of a (smooth or generalized) valuation on a group $G$ and a valuation on a manifold $M$ acted upon by the group. In the case of a transitive group action, we prove that the spaces of smooth and generalized valuations on $M$ are modules over the algebra of compactly supported generalized valuations on $G$ satisfying some technical condition of tameness.

 The case of a vector space acting on itself is studied in detail. We prove explicit formulas in this case and show that the new convolution is an extension of the convolution on smooth translation invariant valuations introduced by J.~Fu and the second named author.
\end{abstract}

\thanks{S.A. was partially supported by ISF grant 1447/12.\\ A. B. was supported by DFG grant BE 2484/5-1.\\ AMS 2010 {\it Mathematics subject
classification}:  53C65, %Integral geometry
22E30%Analysis on real and complex Lie groups
}

\maketitle
%-------------------------------------------------------------------------
\section{Introduction}

\subsection{General background}

A convex valuation on a finite-dimensional vector space $V$ is a map $\mu$ on the space $\mathcal{K}(V)$ of compact convex bodies in $V$ with values in some abelian semi-group such that
\begin{displaymath}
 \mu(K \cap L)+\mu(K \cup L)=\mu(K)+\mu(L)
\end{displaymath}
whenever $K, L, K \cup L \in \mathcal{K}(V)$.

Examples of real valued valuations are any measure on $V$, the Euler characteristic (which equals $1$ on each non-empty $K \in \mathcal{K}(V)$) and the intrinsic volumes in euclidean vector spaces. Examples of valuations with values in other semi-groups are Minkowski valuations \cite{abardia12, abardia_bernig, haberl10, ludwig_2005, schuster10, schuster_wannerer}, tensor valuations \cite{alesker99, bernig_hug, hug_schneider_schuster_a, hug_schneider_schuster_b, mcmullen97}, curvature measures \cite{bernig_fu_solanes,federer59,schneider78, zaehle86} and area measures \cite{schneider75, wannerer_area_measures, wannerer_unitary_module}.  In this paper, we will only consider scalar-valued valuations.

The space of translation invariant and continuous (with respect to Hausdorff metric on $\mathcal{K}(V)$) valuations is denoted by $\Val$. This space carries a surprisingly rich algebraic structure which is moreover very closely related to geometric formulas in Crofton style integral geometry.

Suppose that $V$ is a Euclidean vector space and $G$ a subgroup of $\mathrm{SO}(V)$ acting transitively on the unit sphere in $V$. Then the vector space $\Val^G$ of $G$-invariant elements in $\Val$ is finite-dimensional. Let $\bar G$ denote the group generated by $G$ and translations, endowed with a Haar measure (normalized in an appropriate way). If $\phi_1,\ldots,\phi_N$ is a basis of $\Val^G$, then there exist kinematic formulas of the following type:
\begin{align*}
 \int_{\bar G} \phi_i(K \cap \bar g L) d\bar g & = \sum_{k,l=1}^N c_{k,l}^i \phi_k(K)\phi_l(L), \\
 \int_G \phi_i(K+gL) dg & = \sum_{k,l=1}^N \tilde c_{k,l}^i \phi_k(K) \phi_l(L).
\end{align*}

The first formula is called intersectional kinematic formula, the second one is called additive kinematic formula. In the case where $G=\mathrm{SO}(V)$, a basis of $\Val^G$ is given by the intrinsic volumes, and the coefficients $c, \tilde c$ in the above formulas  may easily be found by plugging in convenient sets $K,L$ (template method). For smaller groups, this method is not strong enough and a more algebraic approach becomes useful.

The space $\Val$ contains a certain dense subspace $\Val^\infty$ of {\it smooth} valuations. The first named author has constructed a product structure on $\Val^\infty$ \cite{alesker04_product}. If $\phi_1,\phi_2 \in \Val^\infty$ are given by
\begin{displaymath}
 \phi_i(K)=\vol_n(K+A_i), \quad i=1,2,
\end{displaymath}
with smooth convex bodies $A_i$ with positive curvatures, then
\begin{displaymath}
 \phi_1 \cdot\phi_2(K)=\vol_{2n}(\Delta(K)+A_1 \times A_2),
\end{displaymath}
where $\Delta:V \to V \times V$ denotes the diagonal embedding.

In the case of a group $G$ as above, $\Val^G$ becomes a subalgebra. J. Fu and the second named author have shown in \cite{bernig_fu06} that the intersectional kinematic formulas and the product structure on $\Val^G$ mutually determine each other. They also introduced a convolution product on $\Val^\infty \otimes \Dens(V^*)$ (with $\Dens(V^*)$ denoting the one-dimensional vector space of densities on $V^*$). If
\begin{displaymath}
 \phi_i(K)=\vol_n(K+A_i) \otimes \vol_n^*, \quad i=1,2,
\end{displaymath}
with smooth convex bodies $A_i$ with positive curvatures, then
\begin{equation} \label{eq_convolution_classic}
 \phi_1 * \phi_2(K)=\vol_n(K+A_1 + A_2) \otimes \vol_n^*.
\end{equation}

The additive kinematic formulas and the convolution product structure on $\Val^G \otimes \Dens(V^*)$ mutually determine each other \cite{bernig_fu06}. Moreover, the first named author constructed a Fourier-type transform $\mathbb{F}: \Val^\infty(V) \to \Val^\infty(V^*) \otimes \Dens(V)$ which intertwines product and convolution and which, in the case of $G$-invariant valuations, gives a direct link between intersectional and additive kinematic formulas \cite{alesker_fourier}. J. Fu coined the name algebraic integral geometry for this area of integral geometry exploiting the links between geometric formulas and algebraic structures \cite{fu_barcelona}.

One success of algebraic integral geometry was the recent determination of intersectional and additive kinematic formulas on hermitian vector spaces under the unitary group $\mathrm{U}(V)$. The first named author constructed a basis of $\Val^{\mathrm{U}(V)}$ \cite{alesker03_un}, J. Fu explicitly computed the algebra structure \cite{fu06} and J. Fu and the second named author translated this algebra structure to a more or less complete set of kinematic formulas \cite{bernig_fu_hig}. Results for other transitive group actions are also available, see \cite{bernig_sun09, bernig_g2,bernig_qig,bernig_solanes, voide}. Wannerer \cite{wannerer_area_measures, wannerer_unitary_module} studied the space of area measures and showed that it is a module over the algebra of smooth and translation invariant valuations, endowed with the convolution product. He determined completely the local additive kinematic formulas for area measures invariant under the group $\mathrm{U}(n)$.

Kinematic formulas do not only exist on affine spaces, but also on some Riemannian manifolds. For instance, on the sphere, one can prove analogues of the classical kinematic formulas which apply to spherically convex bodies moving by the rotation group. On more general Riemannian manifolds, geodesic convexity does not seem to be the right notion. Replacing convexity by a more flexible condition (manifolds with corners), the first named author introduced and studied the space $\mathcal{V}(M)$ of (smooth) valuations on manifolds in a series of papers \cite{alesker_val_man1, alesker_val_man2, alesker_survey07, alesker_val_man4,alesker_val_man3}. An important feature of this space is the existence of a product structure. Very recently, it was shown by J. Fu, G. Solanes and the second named author that the product structure and the kinematic formulas determine each other in the case of isotropic manifolds, i.e. Riemannian manifolds with a subgroup of isometries which acts transitively on the sphere bundle. Using this, they established the full array of kinematic formulas on the complex projective space under its full isometry group \cite{bernig_fu_solanes}.

Another salient feature of the product structure on $\mathcal{V}(M)$ is a Poincar\'e duality which makes it possible to define a space $\mathcal{V}^{-\infty}(M)$ of generalized valuations. For instance, a compact submanifold with corners determines a generalized valuation. In \cite{alesker_bernig}, we constructed a partial product structure on $\mathcal{V}^{-\infty}(M)$. The relation between intersectional kinematic formulas and the product structure is most easily understood in terms of generalized valuations. In fact, the product of two generalized valuations corresponding to compact submanifolds with corners is (under some additional conditions of transversality) the intersection of these submanifolds with corners. Furthermore, generalized valuations appeared naturally in the recent study of valuations which are translation invariant and invariant under the Lorentz group \cite{alesker_faifman}.

In the case of a finite-dimensional vector space, Faifman and the second named author extended the convolution of smooth translation invariant valuations to a partially defined convolution on a certain space of generalized translation invariant valuations. This partial algebra contains McMullen's polytope algebra as a subalgebra (see \cite{bernig_faifman} for the precise statement).

\subsection{Results of the present paper}

As we pointed out in the previous subsection, there are deep connections between algebraic structures on valuations (product, convolution, Fourier transform) and integral-geometric formulas (global and local kinematic formulas on affine and isotropic spaces, additive kinematic formulas). In this paper, we introduce a new type of algebraic structure, which is also called convolution product, which applies to compactly supported valuations on a Lie group, or more generally on a manifold with a transitive Lie group action. We think that it will become a useful tool in the future study of integral geometric formulas on Lie groups and homogeneous spaces.

Let us explain the results in more detail. We consider the action of a Lie group $G$ on a smooth manifold $M$. The convolution of a compactly supported valuation $\mu$ on $G$ with a valuation $\psi$ on $M$ is defined as the push-forward (if it exists) under the multiplication map $a:G \times M \to M$ of the exterior product of $\mu$ and $\psi$ (see the next section for exterior product and push-forward):
\begin{displaymath}
  \mu * \psi :=a_*(\mu \boxtimes \psi).
\end{displaymath}

We will define some special class of generalized valuations, called {\it tame valuations}. The precise definition requires the more technical notion of wave front set of a generalized valuation and is contained in Section \ref{sec_proof_main}.

The space of compactly supported tame valuations on $G$ is denoted by $\mathcal{V}^{-\infty}_{c,t}(G)$. It contains the space of
compactly supported smooth valuations, but is considerably larger. Our main theorem shows that under the assumption that the action of $G$ on $M$ is transitive, the convolution of a tame valuation on $G$ with an arbitrary generalized valuation on $M$ is defined. More precisely, our main theorem is the following.

\begin{MainTheorem} \label{thm_main}
Let $G$ act transitively on $M$. If $\mu$ is a compactly
supported tame generalized valuation on $G$ and $\psi$ a generalized valuation on
$M$, then the convolution $\mu * \psi \in \mathcal{V}_c^{-\infty}(M)$ is
well-defined. Moreover, $\mathcal{V}^{-\infty}_{c,t}(G)$ is an algebra, the spaces
$\mathcal{V}^{-\infty}(M)$, $\mathcal{V}^{\infty}(M)$, $\mathcal{V}^{-\infty}_t(M)$ are modules over
this algebra and the map $(\mu,\psi)\mapsto \mu\ast\psi$ is jointly sequentially continuous in the appropriate topologies.
\end{MainTheorem}

We expect that the assumption on the transitivity of the action can be weakened and that $\mu \ast \psi$ can be defined under more general assumptions.

In the affine case, we can describe the convolution product more explicitly. Moreover, we show that it extends the convolution product of \cite{bernig_fu06}.

\begin{MainTheorem} \label{thm_convolution_rn}
Let $G=V$ be a finite-dimensional vector space acting on itself by addition. Let smooth compactly supported valuations $\phi_i$ on $V$ be given
by
\begin{displaymath}
 \phi_i(K)=\mu_i(K+A_i),\quad i=1,2,
\end{displaymath}
where $\mu_i$ are compactly supported smooth measures on $V$ and $A_i \subset \mathbb{R}^n$ are
smooth compact convex bodies with positive curvature. Then
\begin{equation} \label{eq_defining_equation_flat}
 \phi_1 * \phi_2(K)=\mu_1 * \mu_2(K+A_1+A_2).
\end{equation}
 The convolution product on $\mathcal{V}^\infty_c(V)$ is uniquely characterized by joint sequential continuity, bilinearity and \eqref{eq_defining_equation_flat}. Moreover, there is a natural surjective homomorphism
\begin{displaymath}
 F: (\mathcal{V}^\infty_c(V),*) \to (\Val^\infty(V) \otimes \Dens(V^*),*).
\end{displaymath}
\end{MainTheorem}

%{\red Note the formal similarity with \eqref{eq_convolution_classic}. We do not know whether there is some explicit link between our convolution product on $\mathcal{V}_c(V)$ and the convolution product on $\Val^\infty(V) %\otimes \Dens(V^*)$.}

\subsection{Plan of the paper}

In Section \ref{sec_currents} we collect some well-known facts about currents and wave front sets. Section \ref{sec_genvals} contains the definition and some important properties of generalized valuations. Basic operations such as exterior product, product and push-forward are presented in Section \ref{sec_op_genvals}. Some results on the push-forward appear here for the first time. Section \ref{sec_proof_main} is the heart of the paper. We define the convolution of a generalized valuation on a group and a generalized valuation on a manifold acted upon by the group as the push-forward under the group action of the exterior product. Then we introduce the notion of a {\it tame} generalized valuation and prove existence of the convolution in the case of a transitive group action under the assumption that the valuation on the group is tame. In the final Section \ref{sec_rn} we prove Theorem \ref{thm_convolution_rn}, give examples of tame valuations on $\R^n$ and study in detail the convolution of smooth valuations on $\R$.

\subsubsection*{Thanks}
We wish to thank the Universities of Frankfurt and Tel Aviv for hosting our mutual visits as we worked out this material and Dmitry Faifman and Thomas Wannerer for useful comments.

\section{Currents and wavefronts}
\label{sec_currents}

\subsection{Currents}
Let $X$ be a smooth manifold of dimension $n$. Let $\Omega^k(X)$ denote the space of differential $k$-forms and $\Omega^k_c(X)$ the subspace of compactly supported $k$-forms. We set $\mathcal{D}_k(X):=\Omega^k_c(X)^*$. Elements of this space are called $k$-currents. The boundary of a $k$-current $T$ is the $(k-1)$-current $\partial T$ such that $\langle \partial T,\omega\rangle=\langle T,d\omega\rangle, \omega \in \Omega^{k-1}_c(X)$.  If $\partial T=0$, then $T$ is called a {\it cycle}. A $k$-dimensional closed oriented submanifold $Y \subset X$ (possibly with boundary) defines a current $[[Y]] \in \mathcal{D}_k(X)$ by $\langle [[Y]],\omega\rangle:=\int_Y \omega$. By Stokes' theorem, $\partial [[Y]]=[[\partial Y]]$. If $T \in \mathcal{D}_k(X)$ and $\phi \in \Omega^l(X)$, the current $T \llcorner \phi$ is defined by $\langle T \llcorner \phi,\omega\rangle:=\langle T,\phi \wedge \omega\rangle$.

If $f:X \to Y$ is a smooth map between smooth manifolds $X,Y$ and $T \in \mathcal{D}_k(X)$ such that $f|_{\spt
T}$ is proper, then the push-forward $f_*T \in \mathcal{D}_k(Y)$ is defined by $\langle f_*T,\phi\rangle:=\langle
T,\zeta f^*\phi\rangle$, where $\zeta \in C^\infty_c(X)$ is equal to $1$ in a neighborhood of $\spt T \cap \spt
f^*\phi$.

If $X$ and $Y$ are smooth manifolds, $T \in \mathcal{D}_k(X), S \in \mathcal{D}_l(Y)$, then there is a unique
current $T \boxtimes S \in \mathcal{D}_{k+l}(X \times Y)$ such that $\langle T \boxtimes S,
\pi_1^*\omega \wedge \pi_2^*\phi\rangle=\langle T,\omega\rangle \cdot \langle S,\phi\rangle$, for all $\omega \in
\Omega^k(X), \phi \in \Omega^l(Y)$. Here $\pi_1,\pi_2$ are the projections from $X \times Y$ to $X$ and $Y$
respectively.

\subsection{Wave fronts and operations on currents}

Let $X$ be an oriented manifold of dimension $n$. Let $\p_X:=\p_+(T^*X)$ denote the cosphere bundle, i.e. the $(2n-1)$-dimensional manifold of all tuples $(x,[\xi])$, where $x \in X, \xi \in T_x^*X \setminus \{0\}$ and where the equivalence class $[\xi]$ is with respect to the relation $\xi \equiv \lambda \xi$ for all $\lambda>0$. We denote by $s:\p_X \to \p_X$ the involution $(x,[\xi]) \mapsto (x,[-\xi])$.

The wave front set of a current $T$ on $X$ is a closed conical set $\WF(T) \subset \delbundle{X}$ (cotangent bundle of $X$ with the zero section deleted) describing points and directions of singularities of $T$. In particular, $\WF(T)=\emptyset$ if and only if $T$ is smooth. We refer to \cite{guillemin_sternberg77, hoermander_pde1} for the precise definition of $\WF$. Several operations on currents like pull-back and intersection are only defined under some condition on wave fronts.

If $\Lambda \subset \delbundle{X}$ is a closed conical subset, we let $\bar \Lambda:=\Lambda \cup \underline{0} \subset T^*X$.

\begin{Definition} \label{def_dkgamma}
Let $\Gamma \subset \delbundle{X}$ be a closed conical set. Then  we set
\begin{displaymath}
 \mathcal{D}_{k,\Gamma}(X):=\{T \in \mathcal{D}_k(X): \WF(T) \subset \Gamma\}.
\end{displaymath}
This space has a natural linear locally convex topology (called sometimes H\"ormander topology).
\end{Definition}

\begin{Proposition}[{\cite[Thm. 8.2.3]{hoermander_pde1}}] \label{prop_approx_smooth}
Given $T \in \mathcal{D}_{k,\Gamma}(X)$, there exists a sequence of compactly supported smooth $k$-forms
$\omega_i \in \Omega^{n-k}_c(X)$ such that $[[X]] \llcorner \omega_i \to T$ in $\mathcal{D}_{k,\Gamma}(X)$. In other
words, compactly supported smooth
forms are dense in $\mathcal{D}_{k,\Gamma}(X)$.
\end{Proposition}

\begin{Proposition}[Exterior product] \label{prop_ext_prod}
Let $X,Y$ be smooth manifolds, $\Gamma_1 \subset \delbundle{X}, \Gamma_2 \subset \delbundle{Y}$ closed conical subsets.  Set
\begin{align*}
 \bar \Gamma := \bar \Gamma_1 \times \bar \Gamma_2 \subset T^*(X \times Y).
\end{align*}
Then the exterior product is a jointly sequentially continuous bilinear map
\begin{displaymath}
 \boxtimes:\mathcal{D}_{k,\Gamma_1}(X) \times \mathcal{D}_{l,\Gamma_2}(Y) \to \mathcal{D}_{k+l,\Gamma}(X \times Y).
\end{displaymath}
The joint sequential continuity means explicitly that for sequences $\phi_i\to \phi$
and $\psi_i\to \psi$ we have $\phi_i\boxtimes \psi_i\to \phi\boxtimes\psi$.
\end{Proposition}

\begin{Proposition}[Push-forward] \label{prop_push_forward_current}
Let $f:X \to Y$ be a smooth and proper map between smooth manifolds and $\Gamma \subset \delbundle{X}$ be a closed conical set. Let
\begin{displaymath}
 f_* \Gamma:=\{(y,\eta) \in \delbundle{Y}:\exists x \in f^{-1}(y), (x,df_x^*(\eta)) \in \Gamma\}
\end{displaymath}
 Then the push-forward map $f_*$ of currents is a sequentially continuous map
\begin{displaymath}
 f_*:\mathcal{D}_{k,\Gamma}(X) \to \mathcal{D}_{k,f_*\Gamma}(Y).
\end{displaymath}
\end{Proposition}

\begin{Proposition}[Pull-back] \label{prop_pull_back_current}
Let $f:X \to Y$ be a smooth map between smooth manifolds, where $X$ may have a boundary $N$. Let $\Gamma \subset\delbundle{Y}$ be a closed conical set satisfying the following transversality conditions:
\begin{itemize}
 \item if $x \in X, (f(x),\eta) \in \Gamma$, then $df_x^*(\eta) \neq 0$;
 \item if $x\in N$, $(f(x),\eta)\in \Gamma$, then $(df|_N)^*_x(\eta)\ne 0$.
\end{itemize}
Define
\begin{displaymath}
 f^*\Gamma:=\{((x,df_x^*(\eta)) \in \delbundle{X}: x \in X, (f(x),\eta) \in \Gamma\}.
\end{displaymath}
Then there exists a unique sequentially continuous map
\begin{displaymath}
 f^*:\mathcal{D}_\Gamma(Y) \to \mathcal{D}_{(f^*\bar\Gamma + T^*_NX)\backslash\{\underline{0}\}}(X)
\end{displaymath}
extending the pull-back of smooth forms. It is called pull-back.
\end{Proposition}

The intersection $T_1 \cap T_2$ of two currents $T_1,T_2$ on $X$ is defined as the pull-back (provided it exists) of the current $T_1 \boxtimes T_2$ under the diagonal embedding $\Delta:X \to X \times X$. From Propositions \ref{prop_ext_prod} and \ref{prop_pull_back_current} we obtain the following description of the intersection.

\begin{Proposition}[{\cite[Thm. 8.2.10]{hoermander_pde1}}] \label{prop_intersection_current}
Let $X$ be a smooth manifold of dimension $n$, $\Gamma_1,\Gamma_2 \subset \delbundle{X}$ closed conical sets such that the following transversality condition is
satisfied:
\begin{displaymath}
 \Gamma_1 \cap s \Gamma_2 = \emptyset.
\end{displaymath}
Set
\begin{displaymath}
 \bar \Gamma:=\bar \Gamma_1+\bar \Gamma_2=\left\{(x,\xi_1+\xi_2):(x,\xi_1) \in \bar \Gamma_1, (x,\xi_2) \in
\bar \Gamma_2\right\}.
\end{displaymath}

Then the intersection is a jointly sequentially continuous map
\begin{displaymath}
 \cap: \mathcal{D}_{k_1,\Gamma_1}(X) \times  \mathcal{D}_{k_2,\Gamma_2}(X) \to \mathcal{D}_{k_1+k_2-n,\Gamma}(X).
\end{displaymath}
\end{Proposition}

Let us finally collect some known results which will be important in later sections.

\begin{Proposition} \label{prop_wf_submfld}
 Let $Y \subset X$ be a compact oriented $k$-dimensional submanifold and $[[Y]]$ the $k$-dimensional
current integration against $Y$. Then
\begin{displaymath}
 \overline{\WF}([[Y]])=N_X Y=\{(x,\xi) \in T^*X:\xi|_{T_xY}=0, \quad \forall x \in Y\}.
\end{displaymath}
In particular, with $\delta_x$ being the delta-distribution at $x \in X$, we
get $\overline{\WF}(\delta_x)=T_x^*X$.
\end{Proposition}

\begin{Proposition} \label{prop_wf_restriction_current}
Let $T$ be a $k$-dimensional current on $X$ and $\phi$ an $l$-form on $X$. Then
\begin{displaymath}
 \WF(T \llcorner \phi) \subset \WF(T).
\end{displaymath}
\end{Proposition}

\begin{Proposition}[{\cite[(8.1.11)]{hoermander_pde1}}] \label{prop_wavefront_diffoperator}
Let $P$ be a differential operator on $X$ with infinitely smooth coefficients. Then for any current $T$ on $X$,
\begin{displaymath}
 \WF(PT) \subset \WF(T).
\end{displaymath}
\end{Proposition}

\section{Generalized valuations on manifolds}
\label{sec_genvals}

We assume that $X$ is an oriented $n$-dimensional smooth manifold. The orientability is not necessary and can be easily omitted, but
it simplifies some formulas.
Let $\mathcal{P}(X)$ denote the space of all compact submanifolds with corners. Every $P \in \mathcal{P}(X)$ admits a conormal cycle $N(P)$,
which is a closed Legendrian $(n-1)$-dimensional Lipschitz submanifold of $\p_X$. The orientation of $N(P)$ is fixed in such a way that
\begin{equation} \label{eq_image_down_normal_cycle}
 \pi_* N(P)=\partial P,
\end{equation}
where $\pi:\p_X \to X$ is the projection map.

\begin{Definition}[\cite{alesker_val_man1, alesker_val_man2}]
A smooth valuation $\mu$ on $X$ is a functional $\mu:\mathcal{P}(X) \to \C$ which can be represented in the form
\begin{displaymath}
\mu(P)=\int_{N(P)} \omega+\int_P \phi
\end{displaymath}
with differential forms $\omega \in \Omega^{n-1}(\p_X)$ and $\phi \in \Omega^n(X)$. The space of smooth valuations is denoted by $\mathcal{V}^\infty(X)$ or simply $\mathcal{V}^\infty$ if there is no risk of confusion.
\end{Definition}

The forms $\omega,\phi$ are not uniquely defined by $\mu$. In \cite{bernig_broecker07} it was shown that a pair $(\omega,\phi)$ induces the trivial valuation if and only if
\begin{enumerate}
 \item $D\omega+\pi^*\phi=0$
 \item $\pi_*\omega=0$.
\end{enumerate}
Here $D:\Omega^{n-1}(\p_X) \to \Omega^n(\p_X)$ denotes the Rumin operator \cite{rumin94}. Let us recall its definition. A form
 $\xi$ on a contact manifold is called vertical if for a contact form $\alpha$ one has $\xi \wedge \alpha=0$. Equivalently, $\xi$ is a multiple of $\alpha$: $\xi=\alpha\wedge \eta$ for some form $\eta$. Given $\omega \in \Omega^{n-1}(\p_X)$, there exists a unique vertical form
 $\xi \in \Omega^{n-1}(\p_X)$ such that $d(\omega+\xi)$ is vertical. Then $D\omega:=d(\omega+\xi)$.

The support of a valuation is defined in the obvious way, the space of compactly supported valuations is denoted by $\mathcal{V}_c^\infty$. It admits a natural topology of Fr\'echet space. If $\mu$ is compactly supported, then by \cite[Lemma 2.1.1]{alesker_val_man4}, one can
choose $\omega$ and $\phi$ to be compactly supported as well.

The space $\mathcal{V}^\infty(X)$ carries a natural commutative and associative product. Its construction is involved: first a product on smooth valuations on an affine space is constructed in \cite{alesker04_product} and \cite{alesker_val_man1} (using the solution of P. McMullen's conjecture \cite{alesker_mcullenconj01}). Then it is shown in \cite{alesker_val_man3} that the product can be extended to smooth valuations on an arbitrary manifold $X$ by using local charts. The hard part of this construction is to show that the result is independent of all choices. In \cite{alesker_bernig} we gave another construction which works directly on the level of pairs $(\omega,\phi)$.

The integration functional
\begin{displaymath}
\int: \mathcal{V}^\infty_c(X) \to \mathbb{C}
\end{displaymath}
is, slightly oversimplifying, defined by $\int \mu=\mu(X)$.

The first named author proved the following Poincar\'e duality: the pairing
\begin{align*}
\mathcal{V}^\infty(X) \times \mathcal{V}^\infty_c(X) & \to \mathbb{C},\\
(\mu_1,\mu_2) & \mapsto \int \mu_1 \cdot \mu_2
\end{align*}
is perfect \cite{alesker_val_man4}.

\begin{Definition}[\cite{alesker_val_man4}]
Elements of the space
\begin{displaymath}
\mathcal{V}^{-\infty}(X):=(\mathcal{V}_c^\infty(X))^*
\end{displaymath}
are called generalized valuations on $X$.
\end{Definition}

Note that $\mathcal{V}^\infty(X)$ is a sequentially dense subset of $\mathcal{V}^{-\infty}(X)$ \cite{alesker_bernig}. There is also a canonical injection with dense image
\begin{align*}
\mathcal{P}(X) & \to \mathcal{V}^{-\infty}(X),\\
P & \mapsto [\mu \mapsto \mu(P)].
\end{align*}

\begin{Proposition}[\cite{alesker_bernig}]
 The space $\mathcal{V}^{-\infty}(X)$ of generalized valuations is in one-to-one correspondence with the space of pairs of currents $(T,C)\in \mathcal{D}_{n-1}(\p_X) \times \mathcal{D}_n(X)$ such that $T$ is a Legendrian cycle and $\pi_*T=\partial C$.
\end{Proposition}

More precisely, given $(T,C)$ as above, we may define a generalized valuation $\psi$ by setting
\begin{displaymath}
 \langle \psi,\mu\rangle:=\langle T,\omega\rangle+\langle C,\phi\rangle
\end{displaymath}
for each compactly supported smooth valuation represented by compactly supported forms $\omega,\phi$.

An element $P \in \mathcal{P}(X) \subset \mathcal{V}^{-\infty}(X)$ corresponds to the pair of currents $(N(P),[[P]])$. A smooth valuation, represented by the forms $(\omega,\phi) \in \Omega^{n-1}(\p_X) \times \Omega^n(X)$, corresponds to the pair of currents
\begin{align*}
T & = [[\p_X]] \llcorner s^* (D\omega+\pi^*\phi) \in \mathcal{D}_{n-1}(\p_X),\\
C & = [[X]] \llcorner \pi_*\omega \in \mathcal{D}_n(X).
\end{align*}

\begin{Definition}
 The wave front of a generalized valuation given by a pair $(C,T)$ is defined as the tuple $(\WF(C),\WF(T))$. Given closed conic sets $\Lambda \subset \delbundle{X},\Gamma \subset \delbundle{\p_X}$, we denote by $\mathcal{V}^{-\infty}_{\Lambda,\Gamma}$ the space of generalized valuations $\psi$ such that $\WF(\psi) \subset (\Lambda,\Gamma)$,where the inclusion of pairs of sets is componentwise.
\end{Definition}

\section{Operations on generalized valuations}
\label{sec_op_genvals}
In this section, we collect some known results on exterior product, product and push-forward of generalized
valuations and refer to \cite{alesker_intgeo} for more details. We will also prove a new statement about push-forward (Proposition \ref{prop_push_forward})
which simplifies the construction from \cite[Subsection 3.6.]{alesker_intgeo}  and which will
be used in the sequel.

\subsection{Exterior product}
\label{subsec_extprod}

Let $X_1,X_2$ be smooth manifolds without boundary of dimensions $n_1,n_2$, $X:=X_1 \times X_2$. Let
\begin{align*}
 \mathcal{M}_1 & := \p_{X_1} \times X_2=\{(x_1,x_2,[\xi_1:0]), x_1 \in X_2,x_2 \in X_2, \xi_1 \in T_{x_1}^*X_1 \setminus \{0\}\} \subset \p_X\\
\mathcal{M}_2 & := X_1 \times \p_{X_2}=\{(x_1,x_2,[0:\xi_2]),x_1 \in X_2,x_2 \in X_2, \xi_2 \in T_{x_2}^*X_2 \setminus \{0\}\} \subset \p_X.
\end{align*}

Let $F:\hat \p_X \to \p_X$ be the oriented blow-up of $\p_X$ along $\mathcal{M}:=\mathcal{M}_1 \cup
\mathcal{M}_2$. Let us give an explicit description of $F$.

Let us consider the fiber bundle over $X$ consisting of tuples $(x_1,x_2,[\xi_1:\xi_2],[\xi'_1],[\xi'_2])$
where $[\xi_1:\xi_2] \in \p_+ T^*_{(x_1,x_2)}X$, $[\xi'_1] \in \p_+ T^*_{x_1}X_1$, $[\xi'_2] \in \p_+
T^*_{x_2}X_2$.

Then $\hat \p_X$ is the closure of the set of all such tuples with $\xi_1 \neq 0,\xi_2 \neq 0$ and
$[\xi_1]=[\xi_1'], [\xi_2]=[\xi_2']$. Note that $\hat \p_X$ is a manifold of dimension $2(n_1+n_2)-1$ with
boundary $\mathcal{N}:=\mathcal{N}_1 \cup \mathcal{N}_2$, where
\begin{align*}
 \mathcal{N}_1 & =\{(x_1,x_2,[\xi_1:0],[\xi_1],[\eta_2])\} \\
 \mathcal{N}_2 & =\{(x_1,x_2,[0:\xi_2],[\eta_1],[\xi_2])\} \\
\end{align*}

Then $F$ is given by
\begin{align*}
 F : \hat \p_X & \to \p_X\\
(x_1,x_2,[\xi_1:\xi_2],[\xi_1'],[\xi_2']) & \mapsto (x_1,x_2,[\xi_1:\xi_2]).
\end{align*}
In particular,
\begin{displaymath}
 F(\mathcal{N})=\mathcal{M}=\mathcal{M}_1\cup \mathcal{M}_2 \subset \p_X.
\end{displaymath}

Let $\Phi$ be the map defined by
\begin{align*}
 \Phi: \hat \p_X & \to \p_{X_1} \times \p_{X_2}\\
(x_1,x_2,[\xi_1:\xi_2],[\xi_1'],[\xi_2']) & \mapsto ((x_1,[\xi_1']),(x_2,[\xi_2'])).
\end{align*}

Let
\begin{align*}
 p_1:\p_{X_1} \times X_2 & \to \p_{X_1}\\
p_2:X_1 \times \p_{X_2} & \to \p_{X_2}\\
\tilde p_i:X_1 \times X_2 & \to X_i, i=1,2\\
q_i : \p_{X_1} \times \p_{X_2} & \to \p_{X_i}, i=1,2
\end{align*}
be the natural projections.

Moreover, the inclusions $i_1:\p_{X_1} \times X_2 \to \p_X,i_2:X_1 \times \p_{X_2} \to \p_X$ are defined by
$i_1(x_1,[\xi_1],x_2)=(x_1,x_2,[\xi_1:0])$ and $i_2(x_1,x_2,[\xi_2])=(x_1,x_2,[0:\xi_2])$.

Note that
\begin{displaymath}
 \mathcal{M}=\im i_1 \cup \im i_2
\end{displaymath}

The relevant diagram is given by
\begin{displaymath}
\xymatrix{\p_{X_1}
& \p_{X_1} \times \p_{X_2} \ar@{->>}[l]_-{q_1} \ar@{->>}[r]^-{q_2} & \p_{X_2}\\
& \hat \p_X \ar@{->>}[u]^\Phi \ar@{->>}[d]_F & \\
\p_{X_1} \times X_2 \ar@{->>}[d]_{p_1} \ar@{^{(}->}[r]^-{i_1} & \p_X \ar@{->>}[dd]_{\pi_X} &  X_1 \times \p_{X_2}
\ar@{_{(}->}[l]_-{i_2} \ar@{->>}[d]_{p_2}\\
\p_{X_1} \ar@{->>}[d]_{\pi_{X_1}} & & \p_{X_2} \ar@{->>}[d]_{\pi_{X_2}}\\
X_1 & X_1 \times X_2 \ar@{->>}[l]_-{\tilde p_1} \ar@{->>}[r]^-{\tilde p_2} & X_2}
\end{displaymath}

Now suppose that generalized valuations $\psi_i \in \mathcal{V}^{-\infty}(X_i)$ are given by pairs of
currents $(C_i,T_i)$. Then the generalized valuation $\psi_1 \boxtimes \psi_2 \in \mathcal{V}^{-\infty}(X)$
is defined in \cite[(2.1.12), (2.1.13)]{alesker_intgeo} by the pair of currents
\begin{align}
 C & := C_1 \boxtimes C_2 \label{eq_def_c_extprod}\\
 T & := F_* \Phi^* (T_1 \boxtimes T_2) + (\tilde p_1 \circ \pi_X)^*C_1 \cap (i_{2*}p_2^*T_2)+
(i_{1*}p_1^*T_1) \cap (\tilde p_2 \circ \pi_X)^* C_2 \label{eq_def_t_extprod}.
\end{align}

From this it is easily obtained that
\begin{equation} \label{eq_spt_extprod}
 \spt(\psi_1 \boxtimes \psi_2) \subset \spt(\psi_1) \times \spt(\psi_2).
\end{equation}

It was shown in (\cite[Claim~2.1.9.]{alesker_intgeo}, \cite{alesker_bernig_erratum}) that the exterior product
\begin{displaymath}
 \mathcal{V}^{-\infty}(X_1) \times \mathcal{V}^{-\infty}(X_2) \to
\mathcal{V}^{-\infty}(X_1 \times X_2)
\end{displaymath}
is jointly sequentially continuous. We need the following refinement of this statement.

\begin{Proposition}[Sequential continuity of the exterior product]
Let $\Lambda_i \subset \delbundle{X_i}, \Gamma_i \subset \delbundle{\p_{X_i}}$ be closed conical sets.
Let
\begin{align*}
\bar \Lambda & := \bar \Lambda_1 \times \bar \Lambda_2 \subset T^*(X_1 \times X_2)\\
\bar \Gamma & := F_* (\Phi^* (\bar \Gamma_1 \times \bar \Gamma_2)+T^*_\mathcal{N}\hat \p_X) \\
& \quad \cup [(\tilde p_1 \circ \pi_X)^*\bar \Lambda_1 + i_{2*}p_2^*\bar \Gamma_2]
\cup [i_{1*}p_1^*\bar \Gamma_1+(\tilde p_2 \circ \pi_X)^* \bar \Lambda_2] \subset T^* \p_X.
\end{align*}
Then the exterior product is a jointly sequentially continuous map
\begin{displaymath}
 \mathcal{V}^{-\infty}_{\Lambda_1,\Gamma_1}(X_1) \times \mathcal{V}^{-\infty}_{\Lambda_2,\Gamma_2}(X_2) \to
\mathcal{V}^{-\infty}_{\Lambda,\Gamma}(X_1 \times X_2).
\end{displaymath}
\end{Proposition}

\proof
This follows from \eqref{eq_def_c_extprod}, \eqref{eq_def_t_extprod} and Propositions \ref{prop_ext_prod}, \ref{prop_push_forward_current}, \ref{prop_pull_back_current} and \ref{prop_intersection_current}.
\endproof

\subsection{Product of generalized valuations}
\label{subsec_product_gen}

Let us recall the product of generalized valuations from \cite{alesker_bernig}.

\begin{Theorem}[\cite{alesker_bernig}]
\label{thm_product_gen}
 Let $\Lambda_i \subset \delbundle{X}, \Gamma_i \subset \delbundle{\p_X}$ for $i=1,2$
be closed conical sets. Suppose that the following conditions are satisfied:
\begin{enumerate}
 \item[(a)] $\Lambda_1 \cap s(\Lambda_2)=\emptyset$.
\item[(b)] $\Gamma_1 \cap s(\pi^* \Lambda_2)=\emptyset$, where $\pi=\pi_X:\p_X \to X$ is the natural projection map.
\item[(c)] $\Gamma_2 \cap s(\pi^*\Lambda_1)=\emptyset$.
\item[(d)] If for some $(x,[\xi]) \in \p_X$ we have that $\eta_1 \in \Gamma_1|_{(x,[\xi])}$ and $\eta_2 \in
\Gamma_2|_{(x,[\xi])}$ vanish on the fibers of $\pi$, then $\eta_1 \neq -\eta_2$.
\item[(e)] Let $(x,[\xi],\eta_1) \in \Gamma_1$ with $(x,[\xi]) \in \p_X$ and $\eta_1 \in T^*_{(x,[\xi])}\p_X \setminus
\{0\}$ and $(x,[-\xi],\eta_2) \in \Gamma_2$ with $(x,[-\xi]) \in \p_X$ and $\eta_2 \in T^*_{(x,[-\xi])}\p_X \setminus
\{0\}$.

Let $\tilde \tau:\p_X \times_X \p_X \to \p_X \times \p_X$ be the natural embedding. Set $\zeta:=d\tilde
\tau|_{(x,[\xi],[-\xi])}^*(\eta_1 \times \eta_2) \in T^*_{(x,[\xi],[-\xi])} \p_X \times_X \p_X$. Then
\begin{displaymath}
 d\theta^*(\zeta) \notin T^*_\Delta(\p_X \times_X \p_X)|_{(x,[\xi],[\xi])},
\end{displaymath}
where $\theta(x,[\xi_1],[\xi_2]):=(x,[\xi_1],[-\xi_2])$ and $\Delta$ is the diagonal in $\p_X \times_X \p_X$.
\end{enumerate}
Then the product of smooth valuations extends to a unique jointly sequentially continuous bilinear map
\begin{displaymath}
 \mathcal{V}^{-\infty}_{\Lambda_1,\Gamma_1}(X) \times \mathcal{V}^{-\infty}_{\Lambda_2,\Gamma_2}(X) \to
\mathcal{V}^{-\infty}(X).
\end{displaymath}
\end{Theorem}

\begin{Proposition} \label{prop_continuity_product}
Let $\Lambda \subset \delbundle{X}, \Gamma \subset \delbundle{\p_X}$ be closed conical sets such that
\begin{equation} \label{eq_cond_wave_fronts}
 \pi^*(\Lambda) \subset \Gamma\mbox{ and } \pi^*\pi_*\Gamma\subset \Gamma.
\end{equation}
Then the multiplication map
\begin{displaymath}
\mathcal{V}^\infty(X) \times \mathcal{V}^{-\infty}_{\Lambda,\Gamma}(X) \to
\mathcal{V}^{-\infty}_{\Lambda,\Gamma}(X)
\end{displaymath}
is well-defined and jointly sequentially continuous.
\end{Proposition}

\proof
In \cite{alesker_bernig} we considered the fibre bundle over $X$ consisting of tuples
\begin{displaymath}
 (x,[\xi:\eta],[\xi'],[\eta'],[\zeta]), x \in X, [\xi:\eta] \in \p_+(T_{(x,x)}^*(X \times X)), [\xi'],[\eta'],[\zeta] \in \p_+(T_x^*X).
\end{displaymath}
and defined $\bar \p$ as the closure of the set of all tuples with $\xi,\eta,\xi+\eta \neq 0$ and $[\xi']=[\xi],[\eta']=[\eta],[\xi+\eta]=[\zeta]$.

Then $\bar \p$ is a $(3n-1)$-dimensional manifold whose boundary $\overline{\mathcal{N}}$ consists of the three manifolds
\begin{align*}
 \overline{\mathcal{N}}_0 & := \{(x,[\xi:-\xi],[\xi],[-\xi],[\zeta]), x \in X, \xi,\zeta \in T_x^*X \setminus \{0\}\}\\
 \overline{\mathcal{N}}_1 & := \{(x,[\xi:0],[\xi],[\eta'],[\xi]), x \in X, \xi,\eta' \in T_x^*X \setminus \{0\}\}\\
 \overline{\mathcal{N}}_2 & := \{(x,[0:\eta],[\xi'],[\eta],[\eta]), x \in X, \eta,\xi' \in T_x^*X \setminus \{0\}\}
\end{align*}

The map $\bar \Phi:\bar \p \to \p_X \times_X \p_X$ is defined by
\begin{displaymath}
 \bar \Phi (x,[\xi:\eta],[\xi'],[\eta'],[\zeta])=((x,[\xi']),(x,[\eta']))
\end{displaymath}
and the map $\bar p:\bar \Phi \to \p_X$ by
\begin{displaymath}
 \bar p (x,[\xi:\eta],[\xi'],[\eta'],[\zeta]):=(x,[\zeta]).
\end{displaymath}
Together with the natural projections $q_1,q_2$ we obtain the following diagram.

\begin{displaymath}
 \xymatrix{ & \p_X & \\
 & \bar \p \ar[u]_{\bar p} \ar[d]^{\bar \Phi}& \\
 \p_X & \p_X \times_X \p_X \ar[l]_-{q_1} \ar[r]^-{q_2} & \p_X}
\end{displaymath}

Let $\phi_1 \in \mathcal{V}^\infty(X), \phi_2 \in \mathcal{V}^{-\infty}_{\Lambda,\Gamma}(X)$. Let $(C_i,T_i)$ be
the currents corresponding to $\phi_i$. Then the currents associated to the product $\phi_1 \cdot \phi_2$ are given by

\begin{align*}
C & := C_1 \cap C_2 \in \mathcal{D}_n(X)\\
T & := (-1)^n \bar p_* \bar \Phi^*(q_1^* T_1 \cap q_2^* T_2) + \pi^* C_1 \cap T_2 + T_1 \cap
\pi^* C_2 \in \mathcal{D}_{n-1}(\p_X).
\end{align*}

With $\delta:X \to X \times X$ being the diagonal embedding, we have
\begin{displaymath}
 \overline{\WF}(C) \subset \delta^*(\underline{0} \times \bar \Lambda)=\bar \Lambda.
\end{displaymath}

Hence, since $C_1,T_1$ are smooth we obtain
\begin{displaymath}
\WF(\pi^* C_1 \cap T_2) \subset \Gamma, \WF(T_1 \cap \pi^* C_2) \subset \pi^* \Lambda \subset
\Gamma.
\end{displaymath}

Since $T_1$ is smooth, we have
\begin{displaymath}
 \WF(\bar \Phi^*(q_1^* T_1 \cap q_2^* T_2)) \subset \left((q_2 \circ \bar \Phi)^*\bar\Gamma +
T^*_{\bar{\mathcal{N}}} \bar \p\right) \setminus \underline{0}
\end{displaymath}
provided $\bar{\Phi}^*(q_1^*T_1\cap q_2^*T_2)$ is well defined.
Thus in order to prove the proposition we have to show that
\begin{eqnarray}\label{E:contain1}
(q_2 \circ \bar \Phi)^*\Gamma\cap T^*_{\bar{\mathcal{N}}}\bar \p=\emptyset,\,\\\label{E:contain2}
\bar p_*\left((q_2 \circ \bar \Phi)^*\Gamma +
T^*_{\bar{\mathcal{N}}} \bar \p\right) \subset \Gamma.
\end{eqnarray}

It order to show (\ref{E:contain1})-(\ref{E:contain2}) let us fix a point $\bar\rho\in \bar\p$. Set
$$\rho:=\bar p(\bar\rho)\in \mathbb{P}_X,\, \rho_2:=(q_2\circ \bar \Phi)(\bar\rho)\in \mathbb{P}_X.$$
Next fix
\begin{eqnarray}\label{E:bz}
\bar\zeta\in \left(\bar\Phi^*(q_2^*\Gamma)+T^*_{\bar{\mathcal{N}}}\bar{\p}\right)|_{\bar\rho},\\\label{E:z}
\zeta\in T^*_\rho\mathbb{P}_X \mbox{ such that } d\bar p^*(\zeta)=\bar \zeta.
\end{eqnarray}

In order to check (\ref{E:contain1})-(\ref{E:contain2}) we have to show that
\begin{eqnarray}\label{E:*}
\bar\zeta\ne 0,\\\label{E:**}
\zeta\in \Gamma|_\rho.
\end{eqnarray}

We consider several cases according to which subset of $\bar{\p}$ the point $\bar \rho$ belongs. Below we denote by $\bar\pi\colon \bar{\p}\to X$ the obvious map.

\hfill

\underline{Case 1. Let $\bar \rho\in \bar{\p}\backslash\bar{\mathcal{N}}$.} The condition \eqref{E:*} is satisfied
since $d(q_2\circ \bar \Phi)\colon T_{\bar\rho}\bar{\p}\to T_{\rho_2}\p_X$ is onto. It is easy to see that
\begin{eqnarray}\label{E:c1-1}
\Ker\left(d(q_2\circ \bar\Phi) \colon T_{\bar\rho} \bar{\p}\to T_{\rho_2}\p_X\right)+\Ker\left(d\bar p\colon T_{\bar \rho}\bar{\p}\to T_{\rho}\p_X\right)=\\
\Ker(d\bar\pi \colon T_{\bar\rho}\bar{\p}\to T_{\bar\pi(\bar\rho)}X).
\end{eqnarray}

Next there exists $ \gamma \in \Gamma|_{\rho_2}$ such that
\begin{eqnarray}\label{E:c1-2}
\bar\zeta=d(q_2\circ \bar\Phi)^* \gamma.
\end{eqnarray}

By \eqref{E:c1-1}-\eqref{E:c1-2} and \eqref{E:z} we get that $\bar \zeta\in d\bar\pi^*(T^*X)$. This and \eqref{E:c1-2} imply that
$$\gamma \in  \Gamma|_{\rho_2}\cap d\pi^*(T^*X)|_{\rho_2}.$$
Hence, using the surjectivity of the map $d\bar p|_{\bar \rho}$, one gets $\zeta\in (\pi^*\pi_*\Gamma)|_\rho\subset \Gamma|_\rho$, where the last inclusion is by the assumption of the proposition. Thus the inclusion
\eqref{E:**} is proved.

\hfill

\underline{Case 2. Let $\bar\rho\in \bar{\mathcal{N}_0}$.} We can write
\begin{displaymath}
\bar \zeta=d(q_2\circ \bar\Phi)^*(\gamma)+n=d\bar p^*(\zeta),
\end{displaymath}
where $\gamma\in \Gamma|_\rho,\, n\in T^*_{\bar{\mathcal{N}_0}}\bar\p|_{\bar\rho},\, \zeta\in T^*_\rho\p_X$.

It is easy to see that the restriction
\begin{displaymath}
 (q_2\circ \bar\Phi)|_{\bar{\mathcal{N}_0}}\colon \bar{\mathcal{N}_0}\to \p_X
\end{displaymath}
is a submersion. Hence
\begin{eqnarray}\label{E:c2-1}
\bar\zeta|_{T_{\bar \rho}\bar{\mathcal{N}_0}}=d(q_2\circ \bar\Phi)^*(\gamma)|_{T_{\bar\rho}\bar{\mathcal{N}_0}}\ne 0.
\end{eqnarray}
This implies \eqref{E:*}. Next it is easy to see that
\begin{eqnarray}\label{E:c2-2}
\Ker\left(d(q_2\circ \bar\Phi)|_{T_{\bar\rho}\bar{\mathcal{N}_0}}\colon T_{\bar\rho}\bar{\mathcal{N}_0}\to T_{\rho_2}\p_X\right)+
\Ker \left(d\bar p|_{T_{\bar \rho}\bar{\mathcal{N}_0}}\colon T_{\bar\rho}\bar{\mathcal{N}_0}\to T_{\rho}\p_X\right)=\\\label{E:c2-3}
\Ker(d\bar\pi|_{T_{\bar\rho}\bar{\mathcal{N}_0}}\colon T_{\bar\rho}\bar{\mathcal{N}_0}\to T_{\bar\pi(\bar\rho)}X).
\end{eqnarray}

Clearly \eqref{E:z} and \eqref{E:c2-1}-\eqref{E:c2-3} imply that $\bar\zeta|_{T_{\bar\rho}\bar{\mathcal{N}_0}}\in d\bar\pi^*(T^*X)|_{T_{\bar\rho}\bar{\mathcal{N}_0}}$.
This and \eqref{E:c2-1} imply that
\begin{displaymath}
\gamma \in  \Gamma|_{\rho_2}\cap d\pi^*(T^*X)|_{\rho_2}.
\end{displaymath}
Hence, using the surjectivity of the map $d\bar p|_{\bar \rho}\colon T_{\bar\rho}\bar{\mathcal{N}_0}\to T_\rho \p_X$, one gets
$\zeta\in (\pi^*\pi_*\Gamma)|_\rho\subset \Gamma|_\rho$, where the last inclusion is by the assumption of the proposition. Thus the inclusion
\eqref{E:**} is proved.

\hfill

\underline{Case 3. Let $\bar\rho\in \bar{\mathcal{N}_1}$.} The argument is similar to Case 2. We can write again
\begin{displaymath}
\bar\zeta=d(q_2\circ \bar\Phi)^*(\gamma)+n=d\bar p^*(\zeta),
\end{displaymath}
where $\gamma\in \Gamma|_\rho,\, n\in (T^*_{\bar{\mathcal{N}_1}}\bar{\p})|_{\bar\rho},\, \zeta\in T^*_\rho\p_X$.

It is easy to see that the restriction
\begin{displaymath}
 (q_2\circ \bar\Phi)|_{\bar{\mathcal{N}_1}}\colon \bar{\mathcal{N}_1}\to \p_X
\end{displaymath}
is a submersion. Hence
\begin{eqnarray}\label{E:c3-1}
\bar\zeta|_{T_{\bar\rho}\bar{\mathcal{N}_1}}=d(q_2\circ \bar\Phi)^*(\gamma)|_{T_{\bar\rho}\bar{\mathcal{N}_1}}\ne 0.
\end{eqnarray}
This proves \eqref{E:*}. It is easy to see that
\begin{eqnarray}\label{E:c3-2}
\Ker\left(d(q_2\circ \bar\Phi)|_{T_{\bar\rho}\bar{\mathcal{N}_1}}\colon T_{\bar\rho}\bar{\mathcal{N}_1}\to T_\rho\p_X\right)+
\Ker\left(d\bar p|_{T_{\bar \rho}\bar{\mathcal{N}_1}}\colon T_{\bar \rho}\bar{\mathcal{N}_1}\to T_{\rho_2}\p_X\right)=\\\label{E:c3-3}
\Ker(d\bar\pi|_{T_{\bar\rho}\bar{\mathcal{N}_1}}\colon T_{\bar\rho}\bar{\mathcal{N}_1}\to T_{\bar\pi(\bar\rho)}X).
\end{eqnarray}
Clearly \eqref{E:z} and \eqref{E:c3-1}-\eqref{E:c3-3} imply that
$\bar\zeta|_{T_{\bar\rho}\bar{\mathcal{N}_1}}\in d\bar\pi^*(T^*X)|_{T_{\bar\rho}\bar{\mathcal{N}_1}}$.
This and \eqref{E:c3-1} imply that
\begin{displaymath}
 \gamma \in  \Gamma|_{\rho_2}\cap d\pi^*(T^*X)|_{\rho_2}
\end{displaymath}
and hence $\zeta\in (\pi^*\pi_*\Gamma)|_\rho\subset \Gamma|_\rho$ as in Case 2. Thus the inclusion
\eqref{E:**} is proved.

\hfill

\underline{Case 4. Let $\bar\rho\in \bar{\mathcal{N}_2}.$} We can write again
$$\bar\zeta=d(q_2\circ \bar\Phi)^*(\gamma)+n=d\bar p^*(\zeta),$$
where $\gamma\in \Gamma|_\rho,\, n\in (T^*_{\bar{\mathcal{N}_2}}\bar{\p})|_{\bar\rho},\, \zeta \in T^*_\rho\p_X$.
It is easy to see that the restrictions to $\bar{\mathcal{N}_2}$
\begin{displaymath}
\bar p|_{\bar{\mathcal{N}_2}},\, (q_2\circ \bar\Phi)|_{\bar{\mathcal{N}_2}}\colon \bar{\mathcal{N}_2}\to \p_X
\end{displaymath}
are equal to each other and are submersions. Hence
\begin{eqnarray}\label{E:c4-A}
\bar\zeta|_{T_{\bar\rho}\bar{\mathcal{N}_2}}=d(q_2\circ \bar\Phi)^*(\gamma)|_{T_{\bar\rho}\bar{\mathcal{N}_2}}\ne 0,\\\label{E:c4-B}
\zeta=\gamma.
\end{eqnarray}
Clearly \eqref{E:c4-A} implies \eqref{E:*}, and (\ref{E:c4-B}) implies \eqref{E:**}.

\endproof

\begin{Proposition} \label{prop_density_smooth}
Let $\Lambda \subset \delbundle{X}, \Gamma \subset \delbundle{\p_X}$ be closed cones satisfying \eqref{eq_cond_wave_fronts}.
Then $\mathcal{V}^\infty(X)$ is sequentially dense in $\mathcal{V}^{-\infty}_{\Lambda,\Gamma}(X)$.
\end{Proposition}

\proof
If $X=\R^n$, this was shown in (\cite[Lemma 8.2]{alesker_bernig}). In this case, Condition \eqref{eq_cond_wave_fronts}
is not needed. In the general case, there exist a locally finite open covering $\{U_\alpha\}$ of the manifold $X$ and open sets
$\{\mathcal{O}_\alpha\}$ diffeomorphic to $\mathbb{R}^n$ such that the closure of $U_\alpha$ in $X$ is compact and is contained in $\mathcal{O}_\alpha$.
By \cite[Prop.~ 6.2.1]{alesker_val_man4}, there exists a partition of unity in valuations subordinate to the covering $\{U_\alpha\}$, namely there exist smooth
valuations $\{\phi_\alpha\}$ with $\supp(\phi_\alpha)\subset U_\alpha$ such that
\begin{displaymath}
 \sum_\alpha \phi_\alpha=\chi,
\end{displaymath}
where $\chi$ is the Euler characteristic.

Let $\zeta\in \mathcal{V}^{-\infty}_{\Lambda,\Gamma}(X)$. Then
\begin{displaymath}
\zeta=\sum _\alpha \phi_\alpha\cdot \zeta.
\end{displaymath}
By Proposition \ref{prop_continuity_product} $\phi_\alpha\cdot \zeta\in \mathcal{V}^{-\infty}_{\Lambda,\Gamma}(X)$. Hence we may replace
$\zeta$ by $\phi_\alpha\cdot \zeta$ for some fixed $\alpha$ and assume that $\zeta\in \mathcal{V}^{-\infty}_{\Lambda,\Gamma}(X)$ is a compactly supported valuation
with support contained in an open subset $\mathcal{O}\subset X$ diffeomorphic to $\mathbb{R}^n$.
Then, by the case of $\mathbb{R}^n$, there exists a sequence $\{\zeta_i\}\subset \mathcal{V}^\infty(\mathcal{O})$ which converges to $\zeta|_{\mathcal{O}}$ in $\mathcal{V}^{-\infty}_{\Lambda,\Gamma}(\mathcal{O})$.

Let us choose a smooth compactly supported valuation $\tau$ on $\mathcal{O}$ which is equal to the Euler characteristic $\chi$ in a neighborhood of $\supp(\zeta)$. Then by Proposition \ref{prop_continuity_product}
\begin{displaymath}
\tau \cdot \zeta_i\to \tau\cdot\zeta|_{\mathcal{O}}=\zeta|_{\mathcal{O}}\mbox{ in } \mathcal{V}^{-\infty}_{\Lambda,\Gamma}(\mathcal{O})
\end{displaymath}
and $\supp(\tau\cdot \zeta_i)\subset \supp(\tau)$.

Let $\tilde \zeta_i\in \mathcal{V}^\infty(X)$ be the extension by zero of $\tau\cdot \zeta_i$ from $\mathcal{O}$ to $X$. Then clearly $\tilde\zeta_i\to \zeta$ in $\mathcal{V}^{-\infty}_{\Lambda,\Gamma}(X)$.
\endproof

\subsection{Push forward}

We take up the opportunity to clarify some facts about the push-forward of valuations.
Let $f:X \to Y$ be a smooth proper submersion between smooth manifolds without boundary. Suppose first that $\psi$ is a smooth valuation on $X$. Then the push-forward of $\psi$ under $f$ is defined by
\begin{equation} \label{eq_push_forward_smooth}
 f_*\psi(P):=\psi(f^{-1}P), \quad P \in \mathcal{P}(Y).
\end{equation}

Under some conditions on $\Lambda$ and $\Gamma$ to be specified below, there is a sequentially continuous map
\begin{equation} \label{eq_push_forward_gen}
 f_*:\mathcal{V}^{-\infty}_{\Lambda,\Gamma}(X) \to \mathcal{V}^{-\infty}(Y)
\end{equation}
which extends the push-forward on smooth valuations and which is also called push-forward.

The push-forward $f_* \psi \in \mathcal{V}^{-\infty}(Y)$ of a generalized valuation $\psi \in \mathcal{V}^{-\infty}_{\Lambda,\Gamma}(X)$, if it exists, is defined by the equation
\begin{equation} \label{eq_def_push_forward}
 \int_Y f_*\psi \cdot \phi=\int_X \psi \cdot f^*\phi, \quad \phi \in
\mathcal{V}^\infty(Y).
\end{equation}
Here $f^*:\mathcal{V}^\infty(Y) \to \mathcal{V}^{-\infty}_{\emptyset,T^*_{X \times_{f,\pi_Y}\p_Y}\p_X \setminus \underline{0}}(X)$ denotes
the pull-back of a smooth valuation (see \cite[Section 3.6.]{alesker_intgeo}). The product on
the right hand side is the (partially defined) product of generalized valuations, see Subsection \ref{subsec_product_gen}. The push-forward exists if and only if the right hand side is well-defined for each $\phi$.

The equation \eqref{eq_def_push_forward} implies that
\begin{equation} \label{eq_spt_push_forward}
 \spt f_*\psi \subset f(\spt \phi).
\end{equation}

The map $df^*: X \times_{f,\pi_Y}
\p_Y \to \p_X$ is injective and we will consider $X \times_{f,\pi_Y} \p_Y$ as a subset of $\p_X$.

\begin{Proposition} \label{prop_push_forward}
Let $\psi \in \mathcal{V}^{-\infty}_{\Lambda,\Gamma}(X)$ and $f:X \to Y$. The push-forward $f_*\psi \in
\mathcal{V}^{-\infty}(Y)$ is well-defined provided that
\begin{equation} \label{eq_cond_push_forward}
 \Gamma \cap T^*_{X \times_{f,\pi_Y} \p_Y}\p_X =\emptyset.
\end{equation}
The map
\begin{displaymath}
 f_*:\mathcal{V}^{-\infty}_{\Lambda,\Gamma}(X) \to \mathcal{V}^{-\infty}(Y)
\end{displaymath}
is sequentially continuous.
\end{Proposition}

\proof
We have to show that $\psi \cdot f^*\phi$ is defined for every $\phi \in \mathcal{V}^\infty(Y)$. This amounts to showing that
conditions (a)-(e) in Theorem \ref{thm_product_gen} (with
$(\Lambda_1,\Gamma_1):=(\Lambda,\Gamma)$ and $(\Lambda_2,\Gamma_2):=(\emptyset,T^*_{X \times_{f,\pi_Y}\p_Y}\p_X
\setminus \underline{0})$) are equivalent to \eqref{eq_cond_push_forward}.

Since $\Lambda_2=\emptyset$, conditions (a) and (b) are trivially satisfied.

Since $\pi \circ df^*:X \times_{f,\pi_Y} \p_Y \to X$ is a submersion, it follows that if $\eta \in T^*_{X \times_{f,\pi_Y} \p_Y} \p_X$ vanishes on the fiber of $\p_X$, then $\eta=0$. This in turn implies conditions (c) and (d).

For condition (e), note first that $\Gamma_2$ is closed under $ds^*$.

Let $(x,[\xi],\eta_1) \in \Gamma_1, (x,[-\xi],\eta_2) \in \Gamma_2$ and $\zeta:=d\tilde \tau^*(\eta_1 \times \eta_2)$. Then
\begin{displaymath}
 d\theta^*(\zeta)=d\theta^* \circ d\tilde \tau^*(\eta_1 \times \eta_2)=d\tilde \tau^*(\eta_1 \times \eta_2')
\end{displaymath}
with $\eta_2'=ds^*(\eta_2) \in \Gamma_2|_{x,[\xi]}$.

Now $d\tilde \tau^*(\eta_1 \times \eta_2') \in  T^*_\Delta(\p_X \times_X \p_X)|_{x,[\xi],[\xi]}$ if and only if
$\eta_1=-\eta_2'$. Hence (e) is equivalent to
\begin{displaymath}
 \Gamma \cap T^*_{X \times_{f,\pi_Y}\p_Y}\p_X =\emptyset.
\end{displaymath}
\endproof

We will need a description of the push-forward in terms of currents.

Let $p:X \times_{f,\pi_Y} \p_Y \to \p_Y$ be the canonical projection map. We have a diagram
\begin{displaymath}
 \xymatrix{\p_X & X \times_{f,\pi_Y} \p_Y \ar@{->>}[r]^-p \ar@{(->}[l]_-{df^*}& \p_Y\\
(x,[(d_xf)^* \eta]) & (x,[\eta]), \eta \in T_{f(x)}^*Y \setminus \{0\} \ar@{|->}[l] \ar@{|->}[r] & (f(x),[\eta])}
\end{displaymath}

\begin{Proposition} \label{prop_push_forward2}
Let $f: X\to Y$ be a proper submersion.
Suppose $\Lambda,\Gamma$ satisfy condition \eqref{eq_cond_push_forward}. Let $\psi \in
\mathcal{V}_{\Lambda,\Gamma}^{-\infty}(X)$ be represented by the pair of currents $(C,T)$. Let the
push-forward $f_*\psi$ be represented by the pair of currents $(C',T')$. Then
\begin{equation} \label{eq_push_forward_genval}
 T'=p_* (df^*)^*T.
\end{equation}
\end{Proposition}

Remark: we do not have a simple description of $C'$. However, in some situations $C'$ is uniquely determined by the
relation $\pi_*T'+\partial C'=0$.

\proof
\begin{enumerate}
\item Let us first assume that $\psi$ is smooth, say $\psi$ is represented by a pair $(\omega,\phi)$ of smooth forms.
Then $C=(\pi_X)_*\omega$ and $T=D\omega+\pi_X^*\phi$. By \cite[Prop. 3.2.3]{alesker_intgeo}, the push-forward
$f_*\psi$
is represented by the pair $(f_*\phi,p_*(df^*)^*\omega)$. In the language of currents this means that $f_*\psi$ is
represented by the currents $C':=(\pi_Y)_* p_*(df^*)^*\omega$ and $T':=D(p_*(df^*)^*\omega)+\pi_Y^* f_*\phi$.

We have a commuting diagram
\begin{displaymath}
 \xymatrix{\p_X \ar[rd]_-{\pi_X} & X \times_{f,\pi_Y} \p_Y \ar[r]^-p \ar[l]_-{df^*} \ar[d]^{\pi_1} & \p_Y \ar[d]^{\pi_Y}\\
& X \ar[r]^-f & Y}
\end{displaymath}
where the right hand square is a cartesian square. It follows that
\begin{equation} \label{eq_p*}
 p_*(df^*)^* \pi_X^* \phi=p_* \pi_1^*\phi=\pi_Y^*f_*\phi.
\end{equation}
\item
 We next claim that
\begin{equation} \label{eq_Dp*}
 D(p_*(df^*)^*\omega)=p_*(df^*)^*D\omega
\end{equation}
Since $d$ and $p_*(df^*)^*$ commute, it is enough to show that $p_*(df^*)^*$ maps vertical forms to vertical forms. Let
$(x,[\xi]) \in  X \times_{f,\pi_Y} \p_Y$, i.~e. $x \in X, \xi \in T_{f(x)}^*Y$ with $\xi \neq 0$. Let $\tilde v \in
T_{(x,[\xi])} (X \times_{f,\pi_Y} \p_Y)$.

We have $p(x,[\xi])=(f(x),[\xi]) \in \p_Y$. Then $p_* \tilde v$ is horizontal (i.e. belongs to the contact distribution) if and only if
\begin{displaymath}
\langle \xi,df|_x((\pi_1)_* \tilde v)\rangle=0.
\end{displaymath}
On the other hand, $df^*(x,[\xi])=(x,[df|_x^*(\xi)]) \in \p_X$ and $(df^*)_* \tilde v$ is horizontal if and only if
\begin{displaymath}
\langle df|_x^*(\xi),(\pi_1)_*\tilde v\rangle=0.
\end{displaymath}

Since these conditions are equivalent, we obtain that $(df^*)_*(\tilde v)$ is horizontal in $\p_X$ if (and only if)
$p_*(\tilde v)$ is horizontal
in $\p_Y$. In particular, this is the case if $\tilde v$ is tangent to the fiber of $p$.
From these facts it follows that $p_*(df^*)^* \omega$ is vertical for all vertical forms $\omega$.

From \eqref{eq_p*} and \eqref{eq_Dp*} it follows that
\begin{displaymath}
 p_*(df^*)^*T=p_*(df^*)^*(D\omega+\pi_X^*\phi)=D(p_*(df^*)^*\omega)+\pi_Y^*f_*\phi=T'.
\end{displaymath}

\item Fix $\Lambda,\Gamma$ satisfying condition \eqref{eq_cond_push_forward} and $\psi \in
\mathcal{V}^{-\infty}_{\Lambda,\Gamma}(X)$. We approximate $\psi$ by a sequence of smooth valuations $\psi_i \in \mathcal{V}^\infty(X)$ (see
Proposition \ref{prop_density_smooth}). Let $(C_i,T_i)$ be the corresponding (smooth) currents and $(C_i',T_i')$ the
currents corresponding to $f_*\psi_i$.

Since the map $f_*:\mathcal{V}^{-\infty}_{\Lambda,\Gamma}(X) \to \mathcal{V}^{-\infty}(Y)$ is sequentially continuous
by Proposition \ref{prop_push_forward}, $f_*\psi$ is the limit in $\mathcal{V}^{-\infty}(Y)$ of the sequence $f_*\psi_i$.
In particular, $T_i'$ converges in $\mathcal{D}(Y)$ to $T'$.

On the other hand, by what we have already shown, $T_i'=p_*(df^*)^*T_i$.

The map  $p_*(df^*)^*:\mathcal{D}_\Gamma(X) \to \mathcal{D}(Y)$ is sequentially continuous, compare Section
\ref{sec_currents}. More precisely, the transversality condition in Proposition \ref{prop_pull_back_current} is equivalent
to condition \eqref{eq_cond_push_forward}. We therefore obtain that $T_i'=p_*(df^*)^*T_i  \to p_*(df^*)^* T$ in $\mathcal{D}(Y)$. Finally we get
$T'=p_*(df^*)^* T$.
\end{enumerate}
\endproof

\begin{Corollary} \label{cor_continuity_push_forward}
 In the same situation as in Proposition \ref{prop_push_forward2}, the push-forward map is a continuous map
\begin{displaymath}
 f_*:\mathcal{V}^{-\infty}_{\Lambda,\Gamma}(X) \to \mathcal{V}^{-\infty}_{\Lambda',\Gamma'}(X),
\end{displaymath}
where $\Lambda':=\delbundle{Y}$ and $\Gamma':=p_*(df^*)^*\Gamma \subset \delbundle{\p_Y}$.
\end{Corollary}

\proof
This follows at once from \eqref{eq_push_forward_genval}, taking into account Propositions \ref{prop_push_forward_current} and \ref{prop_pull_back_current}.
\endproof

{\bf Remark:} Suppose that $\Gamma,\Lambda$ satisfy condition \eqref{eq_cond_push_forward}. Set $ \tilde \Gamma:=\pi^*(\Lambda) \cup
\Gamma$. Then $ \tilde \Gamma$ also satisfies condition \eqref{eq_cond_push_forward}. This follows from
the proof of Proposition \ref{prop_push_forward}.

\begin{Proposition}
 Let $f:X \to Y, g:Y \to Z$ be smooth proper submersions between smooth manifolds. Let $\Lambda \subset \delbundle{X}, \Gamma \subset \delbundle{\p_X}$ be closed conical sets satisfying condition \eqref{eq_cond_push_forward}. Let $\Lambda':=\delbundle{Y}$ and $\Gamma':=p_*(df^*)^*\Gamma \subset \delbundle{\p_Y}$ and suppose that $\Gamma'$ satisfies \eqref{eq_cond_push_forward} (with respect to the map $g:Y \to Z$). Then
\begin{displaymath}
 g_* \circ f_* =(g \circ f)_* \text{ on } \mathcal{V}^{-\infty}_{\Lambda,\Gamma}(X).
\end{displaymath}
In particular, both sides of this equation are well-defined maps.
\end{Proposition}

\proof
\begin{enumerate}
 \item Let us first show that $(g \circ f)_*$ is defined on $\mathcal{V}^{-\infty}_{\Lambda,\Gamma}$.

We have the following commutative diagram:
\begin{displaymath}
 \xymatrix{&& \p_Z\\
& X \times_Z \p_Z \ar[d]^q \ar[r]^{f \times \mathrm{id}} \ar[ld]_{d(g \circ f)^*} \ar[ur]^-{p_{XZ}} & Y \times_Z \p_Z \ar[d]^{dg^*} \ar[u]_{p_{YZ}}\\
\p_X & X \times_Y \p_Y \ar[r]_-{p_{XY}} \ar[l]^-{df^*}& \p_Y
}
\end{displaymath}

Here the inclusion $q$ is defined by $q(x,[\eta]):=(x,[dg|_{f(x)}^*\eta])$.

We have to show that $\Gamma \cap T^*_{X \times_Z \p_Z}\p_X = \emptyset$.
We set
\begin{displaymath}
 S:=(df^*)^*\Gamma \subset T^*(X \times_Y \p_Y).
\end{displaymath}
Since $T^*_{X \times_Z \p_Z}\p_X=\ker (d(g\circ f)^*)^*=\ker (df^* \circ q)^*=\ker q^* \circ (df^*)^*$, this amounts to showing that $q^*(\xi) \neq 0$ for all $\xi \in S$.

If this is not the case, then there exists $(x,[\eta]) \in X \times_Z
\p_Z$ and $\xi \in T^*_{q(x,[\eta])}(X \times_Y \p_Y)$ with $\xi \neq 0$, $(q(x,[\eta]),\xi) \in S$ and
$dq|_{x,[\eta]}^* \xi=0$.

The kernel of the map $dp_{XY}|_{q(x,[\eta])}:T_{q(x,[\eta])} (X \times_Y \p_Y) \to T_{p_{XY} \circ q(x,[\eta])}\p_Y$ consists of the vectors of the
form $(v,0)$ with $df|_x(v)=0$. Such vectors are tangent to the image of $q$ and
thus in the kernel of $\xi$. It follows that $\xi=dp_{XY}|_{q(x,[\eta])}^* \xi'$ with $\xi' \in T^*_{p_{XY} \circ
q(x,[\eta])}Y$. This implies that $(f(x),[dg|_{f(x)}^* \eta], \xi') \in {p_{XY}}_*S=\Gamma'$.

Since $\Gamma' \cap T^*_{Y \times_Z \p_Z}\p_Y=\emptyset$ by our assumption, we have $(d|_{(f(x),[\eta])}(dg^*))^*\xi' \neq
0$. Since $f$ is a submersion, we obtain that
\begin{displaymath}
 (d|_{(x,[\eta])}(f \times \mathrm{id}))^* \circ (d|_{(f(x),[\eta])}(dg^*))^*\xi' \neq
0.
\end{displaymath}
This is a contradiction, since the inner square in the above diagram commutes and
\begin{displaymath}
 dq|_{(x,[\eta])}^* \circ
dp_{XY}|_{q(x,[\eta])}^* \xi'=dq|_{(x,[\eta])}^*\xi=0.
\end{displaymath}

It follows that $(g \circ f)_*\psi$ is well-defined for $\psi \in  \mathcal{V}^{-\infty}_{\Lambda,\Gamma}(X)$.
\item
We want to approximate a generalized valuation $\psi \in \mathcal{V}^{-\infty}_{\Lambda,\Gamma}(X)$ by a sequence of smooth ones using Proposition \ref{prop_density_smooth}. For this, we need \eqref{eq_cond_wave_fronts}, which is not satisfied by $\Lambda,\Gamma$ in general. However, setting $\tilde \Gamma:=\Gamma \cup \pi^*(\Lambda)$ and $\tilde \Gamma':=\Gamma' \cup \pi^*(\Lambda')$, the remark following Corollary \ref{cor_continuity_push_forward} implies that the maps
\begin{align*}
f_*:\mathcal{V}^{-\infty}_{\Lambda,\tilde \Gamma}(X) & \to
\mathcal{V}^{-\infty}_{\Lambda',\tilde \Gamma'}(Y)\\
g_*:\mathcal{V}^{-\infty}_{\Lambda',\tilde \Gamma'}(Y) & \to \mathcal{V}^{-\infty}(Z)\\
(g \circ f)_*:\mathcal{V}^{-\infty}_{\Lambda,\tilde \Gamma}(X) & \to \mathcal{V}^{-\infty}(Z)
\end{align*}
are well-defined and continuous.

Let $\psi \in \mathcal{V}^{-\infty}_{\Lambda,\Gamma}(X)$. By Proposition \ref{prop_density_smooth}, we may approximate
$\psi$ in $\mathcal{V}^{-\infty}_{\Lambda,\tilde \Gamma}(X)$ by a sequence of smooth valuations $\phi_i \in
\mathcal{V}^\infty(X)$. Then $g_* \circ f_* \phi_i \to g_*
\circ f_* \psi$ and $(g \circ f)_* \phi_i \to (g \circ f)_* \psi$ in
$\mathcal{V}^{-\infty}(Z)$.

Since $g_* \circ f_* \phi_i=(g \circ f)_* \phi_i$ by \eqref{eq_push_forward_smooth}, it follows that $g_* \circ f_*
\psi=(g \circ f)_* \psi$.
\end{enumerate}
\endproof

%----------------------------------------------------------------------------
\section{Convolution of generalized valuations} \label{sec_proof_main}

\begin{Definition}
 Let $G$ be a Lie group acting on a smooth manifold $M$ by a smooth map $a:G \times M \to M$. Let $\mu$ be a generalized valuation on $G$ and $\psi$ a generalized valuation on $M$. The convolution $\mu * \psi$ is defined as the generalized valuation $a_*(\mu \boxtimes \psi)$, provided that the push-forward exists.
\end{Definition}

In this section, we will show that in the case of a transitive group action and under some extra condition of {\it tameness} of $\mu$, the convolution $\mu*\phi$
is always defined.

First we need some technical lemmas. Let $G$ be a Lie group which acts smoothly and transitively on a smooth manifold
$M$ by $a:G \times M \to M$.

We use the same notation as in Subsection \ref{subsec_extprod}, with $X_1:=G, X_2:=M, X:=G \times M$. Recall that
\begin{align*}
 \mathcal{M}_1 & = \p_G \times M=\{(g,x,[\xi_1:0])| g \in G, x \in M, \xi_1 \in T^*_gG \setminus \{0\}\} \subset \p_{G
\times M}\\
\mathcal{M}_2 & = G \times \p_M=\{(g,x,[0:\xi_2]| g \in G, x \in M, \xi_2 \in T^*_xM \setminus \{0\}\} \subset \p_{G
\times M}
\end{align*}
and $\mathcal{M}=\mathcal{M}_1 \cup \mathcal{M}_2$.

The relevant diagram is
\begin{equation} \label{eq_diagram}
\xymatrix{
\p_G & \p_G \times \p_M \ar@{->>}[l]_-{q_1} \ar@{->>}[r]^-{q_2}  & \p_M & \\
& \hat \p_{G \times M} \ar@{->>}[u]^-\Phi \ar@{->>}[d]_-F & (G \times M) \times_{a,\pi_M} \p_M
\ar@{_{(}->}[ld]_-{da^*} \ar@{->>}[r]^-p \ar[ul]_-r &
\p_M\\
\p_G \times M \ar@{->>}[d]_-{p_1} \ar@{^{(}->}[r]^-{i_1} & \p_{G \times M}  \ar@{->>}[dd]_-{\pi_{G \times M}} &  G \times
\p_M
\ar@{_{(}->}[l]_-{i_2} \ar@{->>}[d]_-{p_2} & \\
\p_G \ar@{->>}[d]_-{\pi_G} & & \p_M \ar@{->>}[d]_-{\pi_M} & \\
G & G \times M  \ar@{->>}[l]_-{\tilde p_1} \ar@{->>}[r]^-{\tilde p_2}  & M &}
\end{equation}
where the map $r$ will be defined below.

\begin{Lemma} \label{lemma_image_astar}
Let $G$ act transitively on $M$.
Let $da^*:(G \times M) \times_{a,\pi_M} \p_M \to \p_{G \times M}$.
Then
\begin{displaymath}
  da^*((G \times M) \times_{(a,\pi_M)} \p_M) \cap \mathcal{M} = \emptyset.
 \end{displaymath}
\end{Lemma}

\proof
Let $(g,x,[\tau]) \in G \times M \times_{a,\pi_M} \p_M$ (i.e. $\tau \in T^*_{gx}M \setminus \{0\}$).
Then
$da^*(g,x,[\tau])=(g,x,[da|_{g,x}^*(\tau)])$. Now suppose that $(g,x,[da|_{g,x}^*(\tau)]) \in \mathcal{M}_1$,
i.e. $[da|_{g,x}^*(\tau)]=[\xi:0]$ for some $\xi \in T_g^*G$.

Let $v \in T_xM$. Then $0=\langle da|_{g,x}^*\tau,(0,v)\rangle=\langle \tau,da|_{g,x}(0,v)\rangle$. However, the map
$T_xM \to
T_{gx}M, v \mapsto da|_{g,x}(0,v)$ is an isomorphism, since $G$ acts by diffeomorphisms. Hence $\tau=0$,
which is a contradiction.

Next, suppose that $(g,x,[da|_{g,x}^*(\tau)]) \in \mathcal{M}_2$,
i.e. $[da|_{g,x}^*(\tau)]=[0:\eta]$ for some $\eta \in T_x^*M$.

Let $u \in T_gG$. Then $0=\langle da|_{g,x}^*\tau,(u,0)\rangle=\langle \tau,da|_{g,x}(u,0)\rangle$. However, the map
$T_gG \to
T_{gx}M, u \mapsto da|_{g,x}(u,0)$ is onto, since $G$ acts transitively. Hence $\tau=0$,
which is a contradiction.
\endproof

By the above lemma, we may define the map
\begin{displaymath}
r:=\Phi \circ F^{-1} \circ da^*:(G \times M) \times_{a,\pi_M} \p_M \to \p_G \times \p_M.
\end{displaymath}

\begin{Lemma} \label{lemma_submersionlemma}
 The map
\begin{displaymath}
 ((\pi_G \circ q_1) \times q_2) \circ r: (G \times M) \times_{a,\pi_M} \p_M \to G
\times \p_M
\end{displaymath}
is a diffeomorphism.
\end{Lemma}

\proof
For $g \in G$ let $\iota_g :M \to G \times M, x \mapsto (g,x)$. Then $a \circ \iota_g=a_g$, i.e.
multiplication by $g$. In particular $da|_{(g,x)} \circ d\iota_g|_x=da_g|_x$ is an isomorphism, from which we deduce that
$d\iota_g|_x^* \circ da|_{(g,x)}^*=da_g|_x^*$ is an isomorphism.

If $(g,x,[\xi]) \in (G \times M) \times_{a,\pi_M} \p_M$, then
\begin{align*}
 ((\pi_G \circ q_1) \times q_2) \circ r(g,x,[\xi]) & =
 ((\pi_G \circ q_1) \times q_2) \circ \Phi \circ F^{-1} \circ da^*(g,x,[\xi]) \\
& =   ((\pi_G \circ q_1) \times
q_2) \circ \Phi \circ F^{-1}(g,x,[da|_{(g,x)}^*\xi])\\
& = (g,x,[d\iota_g|_x^* da|_{g,x}^* \xi])\\
& = (g,x,[da_g|_x^* \xi]).
\end{align*}

From this the statement follows.
\endproof

\begin{Definition}
Let $M$ be a smooth manifold. We define
\begin{displaymath}
 \mathcal{V}^{-\infty}_{t}(M):=\mathcal{V}^{-\infty}_{\delbundle{M},d\pi_M^*(T^*M) \setminus \underline{0}}(M).
\end{displaymath}
Such generalized valuations are called {\it tame}. We also write $\mathcal{V}^{-\infty}_{c,t}(M)$ for
the subspace of compactly supported elements in
$\mathcal{V}^{-\infty}_t(M)$.
\end{Definition}

Note that smooth valuations are tame. Tame valuations appear in a natural way, compare Proposition
\ref{prop_infinitesimal_valuation}.

\begin{Proposition}  \label{prop_existence_convolution}
 Let $G$ act transitively on $M$.
\begin{enumerate}
 \item The convolution product
\begin{equation} \label{eq_convolution_product}
 *: \mathcal{V}_{c,t}^{-\infty}(G) \times \mathcal{V}^{-\infty}(M) \to \mathcal{V}^{-\infty}(M)
\end{equation}
is well-defined and jointly sequentially continuous.
\item Let $\mu \in \mathcal{V}_{c,t}^{-\infty}(G)$ and $\psi \in
\mathcal{V}^{-\infty}(M)$. Let $(C_1,T_1),(C_2,T_2)$ be the currents corresponding to $\mu$ and $\psi$ respectively. Let $(C,T)$ be the currents corresponding to $\mu * \psi$. Then
\begin{equation} \label{eq_conv_intermsof_t}
 T=p_* r^* (T_1 \boxtimes T_2).
\end{equation}
\item $\spt
(\mu*\psi) \subset a(\spt \mu \times \spt \psi)$.
\end{enumerate}
\end{Proposition}

\proof
\begin{enumerate}
 \item Let $\mu \in \mathcal{V}^{-\infty}_{c,t}(G)$ and $\psi \in
\mathcal{V}^{-\infty}(M)$. Let $(C_1,T_1), (C_2,T_2)$ be
the corresponding currents. By Subsection \ref{subsec_extprod}, the exterior product $\mu \boxtimes \psi$
corresponds to currents
\begin{align*}
 C' & := C_1 \boxtimes C_2 \\
 T' & := F_* \Phi^* (T_1 \boxtimes T_2) + (\tilde p_1 \circ \pi_X)^*C_1 \cdot (i_{2*}p_2^*T_2)+
(i_{1*}p_1^*T_1) \cdot (\tilde p_2 \circ \pi_X)^* C_2.
\end{align*}

The convolution of $\mu$ and $\psi$ is defined as the push-forward of $\mu \boxtimes \psi$ (if it exists) under the map $a$. By Proposition \ref{prop_push_forward2}, it corresponds to the currents $(C,T)$ with
\begin{displaymath}
 T=p_* \circ (da^*)^* T'
\end{displaymath}
provided the condition (\ref{eq_cond_push_forward}) is satisfied by $(C',T')$. Let us check it.

Since $\mu \in \mathcal{V}^{-\infty}_{c,t}(G)$, we have
\begin{displaymath}
\overline{\WF}(T_1 \boxtimes T_2)  \subset \overline{\WF}(T_1) \times
\overline{\WF}(T_2) \subset d\pi_G^*(T^*G) \times \overline{\WF}(T_2).
\end{displaymath}

Since $\Phi$ is a submersion, $\Phi^*(T_1 \boxtimes T_2)$ is well-defined.

For a point
$q \in \hat \p \setminus \mathcal{N}$
we get
\begin{align*}
 \overline{\WF}(\Phi^* (T_1 \boxtimes T_2))|_q & \subset d\Phi^*(\overline{\WF}(T_1 \boxtimes
T_2)|_{\Phi(q)}) \\
& \subset d\Phi^*(d\pi_G^*(T^*G) \times \overline{\WF}(T_2)).
\end{align*}

We thus obtain that for $u \in \p_{G \times M} \setminus \mathcal{M}$
\begin{align}
 \overline{\WF}(F_*\Phi^*(T_1 \boxtimes T_2))|_u  & \subset \{d(\Phi \circ
F^{-1})|_u^* (\beta_1,\beta_2) : \nonumber \\
& \quad \quad \beta_1=d\pi_G^*\alpha_1, \alpha_1 \in T^*_{\pi_G \circ q_1 \circ \Phi \circ F^{-1}(u)}G,
\nonumber \\
& \quad \quad  \beta_2 \in
\overline{\WF}(T_2)|_{q_2 \circ \Phi \circ F^{-1}(u)}\}.
\label{eq_first_term_smooth_gen}
\end{align}

Let us set
\begin{displaymath}
A:=\{d(\Phi \circ F^{-1})|_u^* (\beta_1,\beta_2):u \in \p_{G \times M} \setminus \mathcal{M}, \beta_1=d\pi_G^*\alpha_1,
\alpha_1 \in T^*G\} \cup T^
*_{\mathcal{M}}\p_{G \times M}.
\end{displaymath}

From Propositions \ref{prop_ext_prod}, \ref{prop_push_forward_current} and \ref{prop_pull_back_current}, we deduce that the map $(T_1,T_2) \mapsto F_*\Phi^*(T_1 \boxtimes T_2)$ is
jointly sequentially continuous as a map from
$\mathcal{D}_{d\pi_G^*(T^*G)}(\p_G) \times \mathcal{D}(\p_M) \to \mathcal{D}_{A}(\p_{G \times M})$.

Since $i_{2*}p_2^*T_2$ is supported on $\mathcal{M}_2$, we obtain
\begin{equation}
  \overline{\WF}((\tilde p_1 \circ \pi_X)^*C_1 \cap (i_{2*}p_2^*T_2)) \subset T^
*_{\mathcal{M}_2}\p_{G \times M} \subset A.
\label{eq_second_term_smooth_gen}
\end{equation}
The map $(C_1,T_2) \mapsto (\tilde p_1 \circ \pi_X)^*C_1 \cap (i_{2*}p_2^*T_2)$ is a jointly sequentially continuous map
$\mathcal{D}(G) \times \mathcal{D}(\p_M) \to \mathcal{D}_{A}(\p_{G \times M})$.

Similarly, since $i_{1*}p_1^*T_1$ is supported on $\mathcal{M}_1$ we get
\begin{equation}
 \overline{\WF}((i_{1*}p_1^*T_1) \cap (\tilde p_2 \circ \pi_X)^*C_2)|_r \subset T^
*_{\mathcal{M}_1}\p_{G \times M} \subset A.
\label{eq_third_term_smooth_gen}
\end{equation}

The map $(T_1,C_2) \mapsto (i_{1*}p_1^*T_1) \cap (\tilde p_2 \circ \pi_X)^*C_2)$ is a jointly sequentially continuous map
$\mathcal{D}_{d\pi_G^*(T^*G)}(\p_G) \times \mathcal{D}(M) \to \mathcal{D}_{A}(\p_{G \times M})$.

We thus obtain
\begin{align}
 \overline{\WF}(T') & \subset A.
\end{align}

Summarizing this step, we find that the exterior product is a jointly sequentially continuous map
\begin{equation} \label{eq_cont_exterior_product}
\boxtimes: \mathcal{V}^{-\infty}_{c,t}(G) \times \mathcal{V}^{-\infty}(M) \to \mathcal{V}^{-\infty}_{T^*(G \times
M) \setminus \underline{0},A \setminus \underline{0}}(G \times M).
\end{equation}

By Lemma \ref{lemma_image_astar}, the image of $da^*$ is disjoint from $\mathcal{M}$.

Take an element $\lambda \in A \setminus T^*_{\mathcal{M}}(\p_{G \times M})$. Then
$\lambda=d(\Phi \circ
F^{-1})|_u^* (d\pi_G^* \alpha_1,\beta_2)$ with $u \in \p_{G \times M} \setminus \mathcal{M}$,
$\alpha_1 \in T^*_{\pi_G \circ q_1 \circ \Phi \circ F^{-1}(u)}G$ and $\beta_2 \in
\overline{\WF}(T_2)|_{q_2 \circ \Phi \circ F^{-1}(u)}$.

We then obtain

\begin{align*}
 (da^*)^*\lambda & = d(\Phi \circ F^{-1} \circ da^*)^*(d\pi_G^*\alpha_1,\beta_2) \\
& = dr^* \circ d((\pi_G \circ q_1) \times q_2)^*(\alpha_1,\beta_2)\\
& = d(((\pi_G \circ q_1) \times q_2) \circ r)^* (\alpha_1,\beta_2).
\end{align*}
If $\lambda \neq 0$, then $(\alpha_1,\beta_1) \neq (0,0)$ and by Lemma \ref{lemma_submersionlemma} we get
that $(da^*)^* \lambda \neq 0$.

By Proposition \ref{prop_push_forward}, the map
\begin{equation} \label{eq_push_forward2}
 a_*:\mathcal{V}^{-\infty}_{T^*(G \times
M) \setminus \underline{0},A  \setminus \underline{0}}(G \times M) \to \mathcal{V}^{-\infty}(M)
\end{equation}
is defined and sequentially continuous.

From \eqref{eq_cont_exterior_product} and \eqref{eq_push_forward2}, it follows that the convolution product
\eqref{eq_convolution_product} is well-defined and jointly sequentially continuous.
\item The images of $ i_1, i_2$ are contained in $\mathcal{M}$, while the image of $(da^*)^*$ is disjoint from $\mathcal{M}$. Hence the second and third summand in the formula for $T'$ are in the kernel of $(da^*)^*$ and the formula follows by observing that outside $\mathcal{M}$, the map $F$ is a diffeomorphism.
\item The statement about the support of $\mu*\phi$ follows from \eqref{eq_spt_extprod} and
\eqref{eq_spt_push_forward}.
\end{enumerate}
\endproof

\begin{Proposition} \label{prop_submodules}
Let $\mu \in \mathcal{V}^{-\infty}_{c,t}(G)$. If $\psi \in \mathcal{V}^{-\infty}(M)$ is smooth or belongs to
$\mathcal{V}^{-\infty}_t(M)$, then the same holds true for $\mu * \psi$. The maps
\begin{align*}
*: & \mathcal{V}^{-\infty}_{c,t}(G) \times \mathcal{V}^{\infty}(M) \to \mathcal{V}^{\infty}(M) \\
*: & \mathcal{V}^{-\infty}_{c,t}(G) \times \mathcal{V}^{-\infty}_t(M) \to \mathcal{V}^{-\infty}_t(M)
\end{align*}
are jointly sequentially continuous. In particular, if $G$ acts on
itself by multiplication, then $\mathcal{V}^{-\infty}_{c,t}(G)$ is closed under convolution.
\end{Proposition}

\proof
We keep the same notation as in the previous proof. The currents corresponding to $\mu * \psi$ are
$(C,T)$ with $T:=p_*(da^*)^*T'$, which is well-defined by
Proposition \ref{prop_existence_convolution}. Let us compute its wave front.

For $s \in G
\times M \times_{a,\pi_M} \p_M$ we have
\begin{align*}
 \overline{\WF}((da^*)^*T')|_s & \subset \{dr|_s^*(d\pi_G^*\alpha_1,\beta_2):\\
& \quad \quad \alpha_1 \in T^*_{\pi_G \circ q_1 \circ r(s)}G, \\
& \quad \quad \beta_2 \in
\overline{\WF}(T_2)|_{q_2 \circ
r(s)}\}.
\end{align*}

The push-forward $T=p_* (da^*)^*T'$ thus satisfies for $t \in \p_M$
\begin{align}
\overline{\WF}(T))|_t & =  \overline{\WF}(p_* (da^*)^*T'))|_t \\
\quad & \subset \{\eta \in T^*_t\p_M: \exists s \in p^{-1}(t), \exists \alpha_1
\in T^*_{\pi_G \circ q_1 \circ r(s)}G, \nonumber \\
& \quad \quad \exists \beta_2 \in
\overline{\WF}(T_2)|_{q_2 \circ
r(s)} \nonumber \\
& \quad \quad \text{such that } dp|_s^*(\eta) = dr|_s^*(d\pi_G^*\alpha_1,\beta_2)\}. \label{eq_final_result_smooth_gen}
\end{align}

If $\psi$ is smooth, then $T_2$ is smooth and hence $\overline{\WF}(T_2)=\underline{0}$.

Take $\eta$ in the right hand side of \eqref{eq_final_result_smooth_gen}. Then there are $s=(g,x,[\xi]) \in
p^{-1}(t)$ and $\alpha_1 \in
T^*_gG$ with
\begin{displaymath}
 dp|_s^*(\eta) = dr|_s^*(d\pi_G^*\alpha_1,0).
\end{displaymath}

Let $v \in T_t\p_M$. Since $\pi_G \circ q_1 \circ r: G \times M \times_{a,\pi_M}
\p_m \to G$ is the projection on the first factor, while $p$ is projection on the last factor, we may chose a
lift $\tilde v \in T_s(G \times M \times_{a,\pi_m}\p_M)$ such that $dp|_s \tilde v=v$ and $d(\pi_G \circ q_1
\circ r)|_s\tilde v=0$.

Then
\begin{align*}
 \langle \eta,v\rangle & = \langle \eta,dp|_s(\tilde v)\rangle\\
& = \langle dp|_s^*(\eta),\tilde v\rangle\\
& = \langle dr|_s^*(d\pi_G^*\alpha_1,0),\tilde v\rangle\\
& = \langle d(\pi_G \circ q_1 \circ r)^* \alpha_1,\tilde v\rangle\\
& = \langle \alpha_1,d(\pi_G \circ q_1 \circ r)|_s \tilde v\rangle\\
& = 0.
\end{align*}
We conclude that $\eta=0$, which means that $T$ is smooth.

Next, if $\psi \in \mathcal{V}^{-\infty}_t(M)$, then by definition $\overline{\WF}(T_2) \subset
d\pi_M^*(T^*M)$. Let $\eta \in T^*_t\p_M$ be in the right
hand side of \eqref{eq_final_result_smooth_gen}. Then there are $s=(g,x,[\xi]) \in p^{-1}(t)$, $\alpha_1 \in
T^*_gG$, $\alpha_2 \in T^*_xM$ with
\begin{displaymath}
 dp|_s^*(\eta) = dr|_s^*(d\pi_G^*\alpha_1,d\pi_M^*\alpha_2)= d(\Phi \circ F^{-1} \circ
da^*)|_s^*(d\pi_G^*\alpha_1,d\pi_M^*\alpha_2).
\end{displaymath}

Outside $\mathcal{M}$ we have $(\pi_G,\pi_M) \circ \Phi \circ F^{-1}=\pi_{G \times M}$. Using Lemma
\ref{lemma_image_astar} we therefore get
\begin{displaymath}
 dp|_s^*(\eta)=(da^*)^* d\pi_{G \times M}^* (\alpha_1,\alpha_2).
\end{displaymath}

Take a vector $v \in T_t \p_M$ with $d\pi_M(t)=0$. Since $p: G \times M \times_{a,\pi_M} \p_M \to \p_M$ is the
projection on the last factor and $\pi_{G \times M} \circ da^*:  G \times M \times_{a,\pi_M} \p_M \to G \times
M$ is the projection on the first two factors, we may choose a
lift $\tilde v \in T_s (G \times M \times_{a,\pi_M} \p_M)$ such that $dp|_s(\tilde v)=v$ and $d(\pi_{G \times M}
\circ da^*)|_s \tilde v=0$. Then
\begin{align*}
 \langle \eta,v\rangle & = \langle \eta,dp|_s(\tilde v)\rangle\\
& = \langle dp|_s^*\eta,\tilde v\rangle\\
& = \langle[(\alpha_1,\alpha_2),d(\pi_{G \times M}
\circ da^*)|_s \tilde v]\rangle\\
& = 0.
\end{align*}
It follows that $\eta \in d\pi_M^*(T^*M)$. We thus get that $\phi * \psi \in \mathcal{V}^{-\infty}_t(M)$.
\endproof

\begin{Proposition} \label{prop_associativity}
Let $G$ act transitively on $M$.
Let $\mu_1,\mu_2 \in \mathcal{V}^{-\infty}_{c,t}(G)$ and $\phi \in \mathcal{V}^{-\infty}(M)$. Then
\begin{displaymath}
 \mu_1*(\mu_2*\phi)=(\mu_1 * \mu_2) * \phi.
\end{displaymath}
\end{Proposition}

\proof
We need the two actions $m:G \times G \to G$ and $a:G \times M \to M$ of $G$ as well as the map $\tilde a:G \times G \times M \to M, (g_1,g_2,x) \mapsto g_1g_2x$. Correspondingly, we will write $r_G:(G \times G) \times_{m,\pi_G} \p_G \to \p_G \times \p_G$, $r_M:(G \times M) \times_{a,\pi_M} \p_M \to \p_G \times \p_M$, $p_G:(G \times G) \times_{m,\pi_G} \p_G \to \p_G$, $p_M:(G \times M) \times_{a,\pi_M} \p_M \to \p_M$ for the maps in Diagram \eqref{eq_diagram}.

Define maps $r_1: (G \times G \times M) \times_{\tilde a,\pi_M} \p_M \to \p_G \times ((G \times M) \times_{a,\pi_M} \p_M)$ and $r_2:(G \times G \times M) \times_{\tilde a,\pi_M} \p_M \to ((G \times G) \times_{m,\pi_G} \p_G) \times \p_M$ as follows.  Let $(g_1,g_2,x,[\tau]) \in (G \times G \times M) \times_\pi \p_M$, i.e. $\tau \in T_{g_1g_2x}^*M \setminus \{0\}$. If $da|_{g_1,g_2x}^*(\tau)=(\xi_1,\xi_2) \in T_{g_1}^*G \times T_{g_2x}^*M$, then $r_1(g_1,g_2,x,[\tau]):=((g_1,[\xi_1]),(g_2,x,[\xi_2]))$. If $da|_{g_1g_2,x}^*(\tau)=(\xi_1,\xi_2) \in T_{g_1g_2}^* G \times T_x^*M$, then $r_2(g_1,g_2,x,[\tau]):=(g_1,g_2,[\xi_1],x,[\xi_2])$.

Define maps $p_1,p_2:(G \times G \times M) \times_{\tilde a,\pi_M} \p_M \to (G \times M) \times_{a,\pi_M} \p_M$ as follows. Let $(g_1,g_2,x,[\tau]) \in (G \times G \times M) \times_\pi \p_M$, i.e. $\tau \in T_{g_1g_2x}^*M \setminus \{0\}$. Then $p_1(g_1,g_2,x,[\tau]):=(g_1,g_2x,[\tau]), p_2(g_1,g_2,x,[\tau]):=(g_1g_2,x,[\tau])$.

We set
\begin{align*}
 \tilde r_1 & :=(\mathrm{id} \times r_M) \circ r_1,\\
 \tilde r_2 & :=(r_G \times \mathrm{id}) \circ r_2,\\
 \tilde p_1 & :=p_M \circ p_1\\
 \tilde p_2 & :=p_M \circ p_2.
\end{align*}

We then obtain the following two commuting diagrams, where the middle squares are cartesian squares.
\begin{displaymath}
 \xymatrix{(G \times G \times M) \times_{\tilde a,\pi_M} \p_M \ar[d]_-{p_1} \ar[r]^-{r_1} & \p_G \times ((G \times M) \times_{a,\pi_M} \p_M) \ar[d]_-{\mathrm{id} \times p_M} \ar[r]^-{\mathrm{id} \times r_M} & \p_G \times \p_G \times \p_M\\
 (G \times M) \times_{a,\pi_M} \p_M \ar[d]_-{p_M} \ar[r]^-{r_M} &  \p_G \times \p_M &  \\
 \p_M & &  }
\end{displaymath}

\begin{displaymath}
 \xymatrix{(G \times G \times M) \times_{\tilde a,\pi_M} \p_M \ar[d]_-{p_2} \ar[r]^-{r_2} & ((G \times G) \times_{m,\pi_G} \p_G) \times \p_M \ar[d]_-{p_G \times \mathrm{id}} \ar[r]^-{r_G \times \mathrm{id}} & \p_G \times \p_G \times \p_M\\
 (G \times M) \times_{a,\pi_M} \p_M \ar[d]_-{p_M} \ar[r]^-{r_M} &  \p_G \times \p_M &  \\
 \p_M & &  }
\end{displaymath}

Let $(T_1,C_1),(T_2,C_2),(T_3,C_3)$ be the currents representing $\mu_1,\mu_2,\phi$ respectively. Let $(T,C)$ be the currents representing $\mu_1 * (\mu_2 * \phi)$. By \eqref{eq_conv_intermsof_t}
\begin{align*}
 T & =(p_M)_* r_M^* (\mathrm{id} \times p_M)_* (\mathrm{id} \times r_M)^* (T_1 \boxtimes T_2 \boxtimes T_3)\\
 & = (p_M)_* (p_1)_* r_1^* (\mathrm{id} \times r_M)^* (T_1 \boxtimes T_2 \boxtimes T_3)\\
 & = (\tilde p_1)* (\tilde r_1)^* (T_1 \boxtimes T_2 \boxtimes T_3).
\end{align*}

Let $(T',C')$ be the currents representing $(\mu_1 * \mu_2) * \phi$. Then
\begin{align*}
 T' & =(p_M)* r_M^* (p_G \times \mathrm{id})_* (r_G \times \mathrm{id})^* (T_1 \boxtimes T_2 \boxtimes T_3)\\
 & = (p_M)_* (p_2)_* r_2^* (r_G \times \mathrm{id})^* (T_1 \boxtimes T_2 \boxtimes T_3)\\
 & = (\tilde p_2)* (\tilde r_2)^* (T_1 \boxtimes T_2 \boxtimes T_3).
\end{align*}

Since $\tilde p_1=\tilde p_2$ and $\tilde r_1=\tilde r_2$, we deduce that $T=T'$. Hence $\mu_1 * (\mu_2 * \phi)-(\mu_1 * \mu_2) * \phi$ is represented by a pair of currents of the form $(0,\tilde C)$.

Suppose first that $\mu_1,\mu_2,\phi$ are smooth. Then $\tilde C$ is smooth and $\mu_1 * (\mu_2 * \phi)-(\mu_1 * \mu_2) * \phi$ is a multiple of the Euler characteristic valuation on $M$. We want to show that $\mu_1 * (\mu_2 * \phi)=(\mu_1 * \mu_2) * \phi$. By a partition of unity argument, it is enough to prove this under the additional assumption that the supports of the valuations $\mu_1,\mu_2,\phi$ are contained in small open sets $U_1,U_2 \subset G$, $U_3 \subset M$ such that $\tilde a(U_1 \times U_2 \times U_3) \neq M$. Then $\mu_1 * (\mu_2 * \phi)-(\mu_1 * \mu_2) * \phi$ is a multiple of Euler characteristic and supported in a proper subset of $M$, hence it must vanish.

The case of generalized valuations follows by approximation as follows. Let $\mu_1,\mu_2 \in \mathcal{V}^{-\infty}_{c,t}(G), \phi \in \mathcal{V}^{-\infty}(M)$. By Proposition \ref{prop_density_smooth}, there exist sequences $(\mu_1^i),(\mu_2^i)$ in $\mathcal{V}^\infty_c(G)$ and $(\phi^i)$ in $\mathcal{V}^\infty(M)$ converging to $\mu_1,\mu_2,\phi$ in the corresponding topologies. By jointly sequential continuity from Propositions \ref{prop_existence_convolution} and \ref{prop_submodules}, $(\mu_1^i * \mu_2^i) * \phi^i$ converges to $(\mu_1 * \mu_2) * \phi$, while $\mu_1^i * ( \mu_2^i * \phi^i)$ converges to $\mu_1 * (\mu_2 * \phi)$. This finishes the proof.
\endproof

\begin{Corollary} \label{cor_algebra}
The space
\begin{displaymath}
 (\mathcal{V}^{-\infty}_{c,t}(G),*)
\end{displaymath}
is an algebra.
\end{Corollary}

{\bf Example:} For a Lie group $G$ let us denote by $\mathcal{V}^{-\infty}_{\{e\},t}(G)$ the subspace of $\mathcal{V}^{-\infty}_{c,t}(G)$
of generalized valuations supported at identity element $e\in G$. By Proposition \ref{prop_existence_convolution}, $\mathcal{V}^{-\infty}_{\{e\},t}(G)$ is a subalgebra of $\mathcal{V}^{-\infty}_{c,t}(G)$. Furthermore $\mathcal{V}^{-\infty}_{\{e\},t}(G)$ contains the subalgebra of generalized densities (measures) supported at $e$
with respect to the convolution. It is well known that the latter algebra is isomorphic to the universal enveloping algebra $U(\mathfrak g)$ of the Lie algebra $\mathfrak g$ of $G$, i.e.
\begin{displaymath}
 U(\mathfrak g) \subset \mathcal{V}^{-\infty}_{\{e\},t}(G).
\end{displaymath}

\begin{Proposition}
The spaces $\mathcal{V}^{-\infty}(M)$, $\mathcal{V}^{\infty}(M)$ and  $\mathcal{V}^{-\infty}_t(M)$ are
$\mathcal{V}^{-\infty}_{c,t}(G)$-modules.
\end{Proposition}

\proof
This follows from Proposition \ref{prop_associativity} and Corollary \ref{cor_algebra}.
\endproof

\section{Convolution on vector spaces}
\label{sec_rn}
In this section, we study the exterior product as well as the push-forward under the addition map $a:V \times V \to V$ for a finite-dimensional real vector space $V$.

\subsection {Proof of Theorem \ref{thm_convolution_rn}}

A compact convex body $A \subset V$ defines a generalized valuation $\Gamma(A)$ by
\begin{displaymath}
 \langle \Gamma(A),\mu\rangle:=\mu(A), \quad \mu \in \mathcal{V}^\infty_c(V).
\end{displaymath}

\begin{Definition}
 Let $\mu$ be a compactly supported (signed) measure on $V$, $A \subset V$ a compact convex body. We define the generalized valuation $\psi_{\mu,A}$ by
\begin{displaymath}
 \psi_{\mu,A}:=\int_V \Gamma(y-A) d\mu(y) \in \mathcal{V}^{-\infty}(V).
\end{displaymath}
\end{Definition}

Note that, with $\mu$ fixed, the map $A \mapsto \psi_{\mu,A}$ is a continuous map from $\mathcal{K}(V)$ to $\mathcal{V}^{-\infty}(V)$.
\begin{Proposition}
 Let $\mu$ be a smooth compactly supported signed measure, $A$ a smooth compact convex body. Then the image of the smooth valuation $\tau$ with $\tau(K)=\mu(K + A)$ under the injection $\mathcal{V}^\infty(V) \hookrightarrow \mathcal{V}^{-\infty}(V)$ equals $\psi_{\mu,A}$.
\end{Proposition}

\proof
This follows by the definition of the product and the Poincar\'e pairing. If $\nu \in \mathcal{V}^\infty_c(V)$, then
\begin{displaymath}
\langle \tau,\nu \rangle = \int_V \tau \cdot \nu
= \int_V  \nu(y-A)d\mu(y)= \int_V  \langle \Gamma(y-A),\nu\rangle d\mu(y)=\langle \psi_{\mu,A},\nu\rangle.
\end{displaymath}
\endproof

\begin{Proposition} \label{prop_ext_prod_continuous}
Let $V_1,V_2$ be affine spaces. Let $\mu_1,\mu_2$ be compactly supported signed measures and $A_i \in \mathcal{K}(V_i), i=1,2$. Then
\begin{displaymath}
 \psi_{\mu_1,A_1} \boxtimes \psi_{\mu_2,A_2}=\psi_{\mu_1 \boxtimes \mu_2,A_1 \times A_2}.
\end{displaymath}
\end{Proposition}

\proof
Using $\Gamma((y_1-A_1) \times (y_2-A_2))=\Gamma(y_1-A_1) \boxtimes \Gamma(y_2-A_2)$ (\cite[Claim~2.1.11]{alesker_intgeo}), we compute
\begin{align*}
 \psi_{\mu_1 \boxtimes \mu_2,A_1 \times A_2} & = \int_{V \times V} \Gamma((y_1,y_2)-A_1 \times A_2) d(\mu_1 \boxtimes \mu_2)(y_1,y_2)\\
& = \int_V \int_V  \Gamma((y_1-A_1) \times (y_2-A_2)) d\mu_1(y_1) d\mu_2(y_2)\\
& = \int_V \Gamma(y_1-A_1) d\mu_1(y_1)  \boxtimes \int_V  \Gamma((y_2-A_2) d\mu_2(y_2)\\
& = \psi_{\mu_1,A_1} \boxtimes \psi_{\mu_2,A_2}.
\end{align*}
\endproof

\begin{Proposition} \label{prop_push_forward_continuous}
Let $a:V \times V \to V$ be the addition map. Let $\mu$ be a compactly supported smooth signed measure, $A_1,A_2 \in \mathcal{K}(V)$ smooth convex bodies with positive curvature. Then
\begin{displaymath}
 a_* \psi_{\mu,A_1 \times A_2}=\psi_{a_*\mu,A_1+A_2}.
\end{displaymath}
\end{Proposition}

\proof
Let $(B_i)$ be a sequence of smooth convex bodies with positive curvature in $V \times V$ converging to $A:=A_1 \times A_2$. It was shown in \cite[ Prop.~3.6.5.]{alesker_intgeo} that
\begin{displaymath}
a_* \psi_{\mu,B_i}=\psi_{a_*\mu,a(B_i)}.
\end{displaymath}

From $a(B_i) \to A_1+A_2$, it follows that $a_*\psi_{\mu,B_i} = \psi_{a_*\mu,a(B_i)} \to \psi_{a_*\mu,A_1+A_2}$ in
$\mathcal{V}^{-\infty}(V)$.

The set $\mathcal{M}$ of Lemma \ref{lemma_image_astar} is of the form $\mathcal{M}=V \times V \times S \subset \p_{V \times V}$ with
\begin{displaymath}
S:=\{[\xi_1:0]:\xi_1 \in V^* \setminus \{0\}\} \cup \{[0:\xi_2]:\xi_2 \in V^* \setminus \{0\}\} \subset \p_+((V \times
V)^*).
\end{displaymath}

Let $\Lambda \subset T^*(V \times V) \setminus \underline{0}$
be the conic set generated by $\mathcal{M}$. We let
$\Gamma:=\pi_{V \times V}^* \Lambda \subset T^* \p_{V \times V} \setminus \underline{0}$. By Lemma \ref{lemma_image_astar},
\begin{displaymath}
 \Gamma \cap T^*_{(V \times V) \times_{(a,\pi)} \p_V} \p_{V \times V} = \emptyset.
\end{displaymath}

By Proposition \ref{prop_push_forward}, the push-forward map
$a_*:\mathcal{V}^{-\infty}_{\Lambda,\Gamma}(V \times V) \to \mathcal{V}^{-\infty}(V)$ is a sequentially continuous map. If we can chose $B_i$ such that $\psi_{\mu,B_i} \to \psi_{\mu,A_1 \times A_2}$ in
$\mathcal{V}^{-\infty}_{\Lambda,\Gamma}$, then $a_* \psi_{\mu,B_i} \to a_*\psi_{\mu,A_1 \times A_2}$ in the weak topology, and hence
$a_*\psi_{\mu, A_1 \times A_2}=\psi_{a_*\mu,A_1+A_2}$, as claimed.

In order to define $B_i$, we endow $V$ with a Euclidean scalar product and take a sequence of smooth
probability measures $\rho_i$ on $SO(2n)$ whose supports shrink to the unit element. Then we define
$B_i$ by
\begin{displaymath}
 h(B_i,\xi):=\int_{SO(2n)} h(g(A_1 \times A_2),\xi) d\rho_i(g).
\end{displaymath}
$B_i$ is thus the Minkowski integral of the rotated copies of $A_1 \times A_2$. Clearly $B_i$ is a smooth
convex body with positive curvature and $B_i \to A_1 \times A_2$ as $i \to \infty$.

Let $(C,T)$ be the currents corresponding to $\psi_{\mu,A_1 \times A_2}$ and $(C_i,T_i)$ the (smooth) currents corresponding to $\psi_{\mu,B_i}$.
Note that $C_i \to C, T_i \to T$ weakly.

In order to prove convergence of $\psi_{\mu,B_i}$ to $\psi_{\mu,A_1 \times A_2}$ in $\mathcal{V}^{-\infty}_{\Lambda,\Gamma}$, we first
describe $(C_i,T_i)$.

Since $B_i$ is smooth, the normal cycle of $K+B_i$ is given by $(\tau_i)_* N(K)$, where $\tau_i:\p_{V \times V}
\to \p_{V \times V}$ is given by $\tau_i(x_1,x_2,[\xi_1:\xi_2])=((x_1,x_2)+d_{(\xi_1,\xi_2)}
h_{B_i},[\xi_1:\xi_2])$ (\cite[Eq.~(11)]{bernig_quat09}).

Let $\kappa \in \Omega^{2n-1}(V \times V)$ with $d\kappa=\mu$. Then
\begin{displaymath}
\psi_{\mu,B_i}(K)=\int_{K+B_i} \mu=\int_{\partial(K+B_i)} \kappa=\int_{N(K+B_i)} \pi^*\kappa=\int_{N(K)} \tau_i^*
\pi^* \kappa.
\end{displaymath}
We thus have $T_i=D(\tau_i^* \pi^* \kappa)=\tau_i^* \pi^* \mu$.

The support function of $A_1 \times A_2$ is given by
\begin{displaymath}
 h_{A_1 \times A_2}(\xi_1,\xi_2)=h_{A_1}(\xi_1)+h_{A_2}(\xi_2), \quad \xi_1,\xi_2 \in V^*.
\end{displaymath}

It is smooth outside $S$ and
\begin{displaymath}
 h_{B_i}|_{(V \times V)^* \setminus S} \to h_{A_1 \times A_2}|_{(V \times
V)^*
\setminus S}
\end{displaymath}
in $C^\infty((V \times V)^* \setminus S)$.
Therefore $\tau_i$ converges in $C^\infty(\p_{V \times V} \setminus \mathcal{M},\p_{V \times V} \setminus \mathcal{M})$. Hence the
currents $T_i$ converge smoothly to $T$ outside $\Gamma$.

The current $C_i$ is given by the smooth function $1_{B_i} * \mu \in C_c^\infty(V)$, which converges in $C_c^\infty(V)=\mathcal{D}_{n,\emptyset}(V)$ to the function $1_{A_1 \times A_2} * \mu$.

It
follows that $\psi_{\mu,B_i} \to \psi_{\mu,A_1 \times A_2}$ in $\mathcal{V}^{-\infty}_{\Lambda,\Gamma}(V \times V)$.
\endproof

\proof[Proof of Theorem \ref{thm_convolution_rn}]
 Equation \eqref{eq_defining_equation_flat} is immediate from Propositions \ref{prop_ext_prod_continuous} and \ref{prop_push_forward_continuous}. Valuations of the type $K \mapsto \mu(K+A)$ with $\mu$ a compactly supported smooth measure and $A$ a compact convex body with smooth boundary and positive curvature span a dense subset in $\mathcal{V}^\infty_c(V)$ (see \cite[Corollary 3.1.7.]{alesker_val_man1}), which implies the uniqueness statement.

The map $F:\mathcal{V}^\infty_c(V) \to \Val^\infty(V) \otimes \Dens(V^*)$ was introduced in \cite{bernig_faifman}:
\begin{displaymath}
 F(\phi)=\int_V \phi(\cdot+x) d\vol(x) \otimes \vol^*.
\end{displaymath}
If $\phi(K)=\mu(K+A)$ with a smooth measure $\mu$ and a smooth compact convex body $A$ with positive curvature, then one easily checks that $F\phi(K)=\mu(V) \vol(K+A) \otimes \vol^*$.

Let us show that $F$ is surjective. Start with a smooth translation invariant valuation $\mu$. Write $\mu$ in terms of forms $(\omega,\phi)$. Take a smooth compactly supported function $f$ on $V$. Then the valuation represented by forms ($\pi^*f \omega, f\phi)$ has compact support, and it is mapped to $\mu$ times $\int f(x) d\vol(x) \otimes \vol^*$. Hence it is enough to choose f such that its integral is non-zero.

Let now $\phi_i(K)=\mu_i(K+A_i), i=1,2$ as above. Then
\begin{align*}
 F(\phi_1 * \phi_2) & = F((\mu_1 * \mu_2)(\cdot+A_1+A_2))\\
 & = (\mu_1 * \mu_2)(V) \vol(\cdot+A_1+A_2) \otimes \vol^*\\
 & = [\mu_1(V) \vol(\cdot + A_1) \otimes \vol^*] * [\mu_2(V) \vol(\cdot + A_2) \otimes \vol^*] \\
 & = F(\phi_1) * F(\phi_2).
\end{align*}
This shows that $F$ is indeed a homomorphism of algebras and finishes the proof.

\endproof

\subsection{Tame valuations on a vector space}

\begin{Proposition} \label{prop_infinitesimal_valuation}
 Let $K \subset \R^n$ be a smooth compact convex body. Then the
generalized valuation $\tau_K$ given by
\begin{displaymath}
\langle \tau_K,\mu\rangle:=\left.\frac{d}{dt}\right|_{t=0} \mu(tK),\quad  \mu \in \mathcal{V}_c^\infty(\R^n)
\end{displaymath}
is tame.
\end{Proposition}

\proof
Let us assume that $n \geq 2$. Let us do computations in coordinates, identifying $(\R^n)^*$ and $\R^n$.

Let $h_K:\R^n \to \R$ be the support function of $K$. By assumption, $h$ is smooth outside the origin. Let
$g=(g_1,\ldots,g_n):=\grad h|_{S^{n-1}}:S^{n-1} \to \R^n$ and set
\begin{displaymath}
 G_t:=(tg,\mathrm{id}):S^{n-1} \to \R^n \times S^{n-1}.
\end{displaymath}
Then
\begin{displaymath}
 N(tK)=(G_t)_*S^{n-1}.
\end{displaymath}
Let $\mu \in \mathcal{V}_c^\infty(\R^n)$ be represented by the pair $(\omega,\phi)$ of compactly supported forms. Then

\begin{align*}
 \mu(tK) & =\int_{N(tK)}\omega+\int_{tK}\phi\\
& = \int_{(G_t)_*(S^{n-1})} \omega+\int_{tK}\phi\\
& = \int_{S^{n-1}} G_t^*\omega+\int_{tK}\phi\\
\end{align*}

Let us decompose
\begin{displaymath}
 \omega=\sum_{I,j} f_I dx_I \wedge \kappa_I,
\end{displaymath}
where $I$ ranges over all multi-indices of order $\leq n-1$ in $\{1,\ldots,n\}$, $\kappa_I \in \Omega^{n-1-\#I}(S^{n-1})$ and $f_I \in C_c^\infty(\R^n \times S^{n-1})$.

Then
\begin{align*}
 \mu(tK) & = \int_{S^{n-1}} f_I(tg,\cdot) \bigwedge_{i \in I} d(tg_i) \wedge \kappa_I+\int_{tK}\phi\\
& = \int_{S^{n-1}} t^{\#I}f_I(tg,\cdot) \bigwedge_{i \in I} dg_i \wedge \kappa_I+\int_{tK}\phi.
\end{align*}

Taking the derivative at $t=0$ gives us
\begin{align}
 \langle \tau_K,\mu\rangle & =\left.\frac{d}{dt}\right|_{t=0} \mu(K) \nonumber\\
& =\sum_{j=1}^n \int_{S^{n-1}} \frac{\partial
f_\emptyset}{\partial x_j}(0,\cdot)g_j \kappa_\emptyset+\sum_{i=1}^n \int_{S^{n-1}}
f_{\{i\}}(0,\cdot)dg_{i} \wedge \kappa_{\{i\}}. \label{eq_tau_k}
\end{align}

If $\rho$ is any smooth $(n-1)$-form on the unit sphere and $j \in \{1,\ldots,n\}$, then the
distribution
\begin{displaymath}
 T(f):=\int_{S^{n-1}} \frac{\partial f(0,\cdot)}{\partial x_j} \rho=\int_{\{0\} \times S^{n-1}} \frac{\partial
f}{\partial x_j} \pi_2^*\rho, \quad f \in C_c^\infty(\R^n \times S^{n-1})
\end{displaymath}
equals $\pm \frac{\partial}{\partial x_j}([[\{0\} \times S^{n-1}]] \llcorner \pi_2^*\rho)$. By Propositions
\ref{prop_wf_submfld}, \ref{prop_wf_restriction_current} and \ref{prop_wavefront_diffoperator}
\begin{align*}
 \overline{\WF}(T) & \subset \overline{\WF}([[\{0\} \times S^{n-1}]] \llcorner \pi_2^* \rho) \\
& \subset
\overline{\WF}\left([[\{0\} \times S^{n-1}]]\right)\\
& =N_{\{0\} \times S^{n-1}} (\R^n \times S^{n-1})\\
& =d\pi|_{(0,v)}^*(T^*_0\R^n).
\end{align*}
For similar reasons, the wave front coming from the second term in \eqref{eq_tau_k} is contained in
$d\pi^*(T^*_0\R^n)$. We deduce that $\tau_K \in \mathcal{V}_{t}^{-\infty}(\R^n)$.
\endproof

The smoothness assumption in the previous proposition can not be dropped, as is shown by the following example.

\begin{Proposition}
 Let $Q=[-1,1] \times [-1,1] \subset \R^2=\C$. Then $\tau_Q$ is not tame.
\end{Proposition}

\proof
Let $\mu(K):=\int_{N(K)} \omega$, where $\omega:=a_1 dx_1+a_2 dx_2$ with smooth functions $a_1,a_2$ on $\C
\times S^1$. Then
\begin{displaymath}
\mu(t Q)=\int_{N(t Q)}\omega=\int_{-t}^t a_1(s,t,-i)ds+\int_{-t}^t
a_2(t,s,1)ds-\int_{-t}^ta_1(s,t,i)ds-\int_{-t}^ta_2(t,s,-1)ds
\end{displaymath}
and therefore
\begin{align*}
 \left.\frac{d}{dt}\right|_{t=0}
\mu(tK) & =2a_1(0,-i)+2a_2(0,1)-2a_1(0,i)-2a_2(0,-1)\\
& =2(\delta_{(0,-i)}-\delta_{(0,i)})a_1+2(\delta_{(0,1)}-\delta_{(0,
-1)})a_2.
\end{align*}
By Proposition \ref{prop_wf_submfld},
\begin{displaymath}
 T^*_{(0,1)}\p_{\R^2} \cup T^*_{(0,-1)}\p_{\R^2} \cup T^*_{(0,i)}\p_{\R^2} \cup
T^*_{(0,-i)}\p_{\R^2} \subset \overline{\WF}(\tau_Q),
\end{displaymath}
which implies that $\tau_Q \not\in \mathcal{V}^{-\infty}_t(\R^n)$.
\endproof

\subsection{Convolution on $\R$}

Now let us work out in more detail the special case $n=1$, i.e. $G=\R$ acting on itself by multiplication.

It is easy to check that each $\phi \in \mathcal{V}^\infty(\R)$ has the form
\begin{displaymath}
 \phi([a,b])=g(b)-f(a)
\end{displaymath}
with $f,g \in C^\infty(\R)$. These functions are unique up to an addition of the same constant. The support of $\phi$ is compact if
and only if there are $c_1,c_2 \in \R$ with $f(x)=g(x)=c_1$ for all $x$ near $-\infty$ and $f(x)=g(x)=c_2$ for
all $x$ near $+\infty$. In this case, $\int_\R \phi=c_2-c_1$. Note that $f,g$ can not be chosen to be
compactly supported in general. However, we can achieve that $c_1=0$ and we will do this in the following.

In this case, the convolutions $f_1 * f_2, g_1 * g_2$ are defined in the usual way, e.g.
\begin{displaymath}
 f_1 * f_2(x):=\int_\R f_1(y)f_2(x-y)dy.
\end{displaymath}
Note that $f_1 * f_2, g_1 * g_2$ vanish near $-\infty$.

\begin{Proposition}
Let $\phi_i \in \mathcal{V}^\infty_c(\R), i=1,2$ correspond to the pair $(f_i,g_i)$ with $f_i,g_i \in
C^\infty(\R), f_i=g_i=0$ near
$-\infty$ and $f_i=g_i=const$ near $\infty$. Then $\phi_1 * \phi_2$ corresponds to the pair
\begin{displaymath}
 ((f_1 * f_2)',(g_1*g_2)').
\end{displaymath}
\end{Proposition}

\proof
Suppose first that $\int_\R \phi_1=\int_\R \phi_2=0$. Then $f_i,g_i$ are compactly supported functions.
Define compactly supported measures $\gamma_i:=g_i(x)dx, \eta_i:=f_i(x)dx$. By our assumption,
$\gamma_i,\eta_i$ are compactly supported, hence the convolutions $\gamma_1 * \gamma_2, \gamma_1 * \eta_2, \eta_1 * \gamma_2, \eta_1 * \eta_2$ are
well-defined and admit $g_1 * g_2, g_1 * f_2, f_1 * g_2, f_1 * f_2$ as densities.

Let $\epsilon,\epsilon_1,\epsilon_2>0$.

We have
\begin{align*}
 \phi_i([a,b]) & = g_i(b)-f_i(a) \\
& = \left.\frac{d}{d\epsilon}\right|_{\epsilon=0} \int_a^{b+\epsilon}
g_i(x)dx-\left.\frac{d}{d\epsilon}\right|_{\epsilon=0} \int_{a-\epsilon}^b f_i(x)dx \\
& =\left.\frac{d}{d\epsilon}\right|_{\epsilon=0}
\gamma_i([a,b]+[0,\epsilon])-\left.\frac{d}{d\epsilon}\right|_{\epsilon=0} \eta_i([a,b]+[-\epsilon,0]).
\end{align*}

By Theorem \ref{thm_convolution_rn} it follows that
\begin{align*}
 \phi_1*\phi_2([a,b]) &
=\left.\frac{d}{d\epsilon_1}\right|_{\epsilon_1=0}\left.\frac{d}{d\epsilon_2}\right|_{\epsilon_2=0}
\gamma_1*\gamma_2([a,b+\epsilon_1+\epsilon_2])\\
& \quad -\left.\frac{d}{d\epsilon_1}\right|_{\epsilon_1=0}\left.\frac{d}{d\epsilon_2}\right|_{\epsilon_2=0}
\gamma_1*\eta_2([a-\epsilon_2,b+\epsilon_1])\\
& \quad -\left.\frac{d}{d\epsilon_1}\right|_{\epsilon_1=0}\left.\frac{d}{d\epsilon_2}\right|_{\epsilon_2=0}
\eta_1*\gamma_2([a-\epsilon_1,b+\epsilon_2])\\
& \quad +\left.\frac{d}{d\epsilon_1}\right|_{\epsilon_1=0}\left.\frac{d}{d\epsilon_2}\right|_{\epsilon_2=0}
\eta_1*\eta_2([a-\epsilon_1-\epsilon_2,b])\\
&
=\left.\frac{d}{d\epsilon_1}\right|_{\epsilon_1=0}g_1*g_2(b+\epsilon_1)+\left.\frac{d}{d\epsilon_1}\right|_{
\epsilon_1=0 }
f_1*f_2(a-\epsilon_1)\\
& = (g_1*g_2)'(b)-(f_1*f_2)'(a).
\end{align*}

Next let $\phi_1,\phi_2$ be arbitrary. We want to show that
\begin{displaymath}
 \phi_1*\phi_2([a,b])=(g_1*g_2)'(b)-(f_1*f_2)'(a)
\end{displaymath}
for all $a<b$.

Let $N_0$ be such that $\spt \phi_1,\spt \phi_2 \subset [N_0,\infty)$ and fix $N_1>N_0$
such that $N_0+N_1>b$.

Let $\tilde g_i,\tilde f_i$ be {\it compactly supported} smooth functions which agree with $g_i,f_i$ on
$(-\infty,N_1]$ and define the corresponding valuations $\tilde \phi_i \in \mathcal{V}^\infty_c(\R)$. Note
that $\int_\R \tilde \phi_i=0$.

Then, using Proposition \ref{prop_existence_convolution},
\begin{align*}
\spt(\phi_1*\phi_2-\tilde \phi_1*\tilde \phi_2) & =\spt((\phi_1-\tilde \phi_1)*\phi_2+\tilde
\phi_1*(\phi_2-\tilde \phi_2))\\
& \subset  (\spt(\phi_1-\tilde \phi_1)+\spt \phi_2) \cup (\spt \tilde
\phi_1+\spt (\phi_2-\tilde \phi_2))\\
& \subset [N_0+N_1,\infty).
\end{align*}
Since $b<N_0+N_1$, we have in particular $\phi_1*\phi_2([a,b])=\tilde \phi_1*\tilde \phi_2([a,b])$. In a
similar way, one gets $(f_1*f_2)'(a)=(\tilde f_1*\tilde f_2)'(a)$ and $(g_1*g_2)'(b)= (\tilde g_1*\tilde
g_2)'(b)$.

By the case which we have already treated, it follows that
\begin{align*}
\phi_1*\phi_2([a,b]) & =\tilde \phi_1*\tilde \phi_2([a,b])\\
& = (\tilde g_1*\tilde g_2)'(b)-(\tilde f_1*\tilde f_2)'(a)\\
& = (g_1*g_2)'(b)-(f_1*f_2)'(a).
\end{align*}
\endproof

\def\cprime{$'$}

%\bibliographystyle{plain}
%\bibliography{../biblio}

\end{document}